%% file: main.tex
\documentclass[12pt,twoside,a4paper]{article}
\usepackage{article}
\usepackage{notation}

\newcommand{\papertitle}{The construction of Dirac operators\protect\\ on Orientifolds}
\newcommand{\paperhead}{The construction of Dirac operators on Orientifolds}
\newcommand{\paperauthor}{Simon Kitson}

\begin{document}
\color{textcolor}

\pagestyle{fancy}
\fancyhead{}
\fancyhead[EC]{\fontsize{9}{12}\selectfont{\MakeUppercase{\paperauthor}}}
\fancyhead[OC]{\fontsize{9}{12}\selectfont{\MakeUppercase{\paperhead}}}
\fancyhead[EL]{\fontsize{9}{12}\selectfont{\thepage}}
\fancyhead[OR]{\fontsize{9}{12}\selectfont{\thepage}}
\fancyfoot{}
\renewcommand{\headrulewidth}{0pt}
\renewcommand{\footrulewidth}{0pt}

\title{\normalsize\bf\MakeUppercase{\papertitle}}
\author{\fontsize{10}{12}\selectfont\MakeUppercase{\paperauthor}\thanks{I would like to thank the Mathematical Sciences Institute of the Australian National University for the postdoctoral fellowship which supported this research.}}

\date{\vspace{-1cm}}
\maketitle

\begin{center}
    \begin{minipage}{0.875\textwidth}
    {\scshape\noindent Abstract.}{\ \small\input{abstract}}
    \end{minipage}
\end{center}

\begin{center}
\begin{minipage}{0.875\textwidth}
\tableofcontents 
\end{minipage}
\end{center}

\input{body}

\bibliographystyle{plain} 
\bibliography{/home/simon/Library/Library-Contents}


\end{document}

%% file: abstract.tex
%
%
Motivated by Wigner's theorem,
    a canonical construction is described
    that produces an Atiyah-Singer Dirac operator
    with both 
    unitary and anti-unitary symmetries.
    This Dirac operator includes
    the Dirac operator for $KR$-theory 
    as a special case, 
    filling a long-standing gap 
    in the literature.
In order to make the construction,
    orientifold $\Spinc$-structures
    are defined and classified
    using
    semi-equivariant Dixmier-Douady theory,
    and analogues 
    of several standard theorems
    on the existence 
    of $\Spinc$-structures are proved.

%% file: body.tex
%
%

\section{Introduction}

This paper 
    uses new results on
    semi-equivariant Dixmier-Douady theory
    \cite{Kitson.2020.A-Semi-equivariant-Dixmier-Douady-Invariant}
    to determine the orientiation
    conditions that allow the construction
    of Atiyah-Singer Dirac operators
    with both linear and anti-linear symmetries.
    The construction will be described 
    in detail.
The existence of Dirac operators 
    with linear/anti-linear symmetries 
    is basic to the compatibility
    between index theory and Wigner's Theorem. 
    It also fills the gap 
    in the literature
    regarding the existence 
    of Dirac operators for $KR$-theory.
    In addition,
    it seems likely that
    such operators 
    have important applications in
    orientifold string theories
    and
    condensed matter physics.

    In the present context, 
    the term \emph{orientifold} will refer 
    to a manifold 
    equipped with an action of a group $\Gam$
    which, in turn, is equipped 
    with a homomorphism 
        $\eps: \Gam \rightarrow \Z_2$.
    This small amount of extra structure 
    is used to define 
    unitary/anti-unitary actions of $\Gam$ 
    on complex vector bundles 
    over the orientifold. 
    An element $\gam \in \Gam$ 
    acts via 
    a unitary map if 
    $\gam \in \Gam^+ := \ker(\eps)$, 
    or 
    an anti-unitary map if 
    $\gam \in \Gam^- := \Gam \setminus \Gam^+$.
    These vector bundles will be described as \emph{orientifold bundles}.
    Note that the set of orientifold bundles over an orientifold
    depends on the embedding $\Gam^+ \hookrightarrow \Gam$.
    More generally, 
    the term \emph{orientifold} will be used as an adjective to describe
    objects carrying, or compatible with, unitary/anti-unitary actions.
    For example, 
    the Dirac operator
    mentioned above
    acts between orientifold bundles
    in an equivariant manner
    and will be described as the
    \emph{orientifold Dirac operator}.

The construction 
    of the orientifold Dirac operator 
    depends on an understanding 
    of the global topology 
    of complex vector bundles 
    with anti-unitary symmetries.
The main obstacle 
    to understanding the conditions 
    under which 
    an orientifold Dirac operator exists
    is the failure of 
    equivariant transition cocycles 
    and cohomology
    to accomodate anti-linear symmetries.
    This obstacle was overcome in 
    \cite{Kitson.2020.A-Semi-equivariant-Dixmier-Douady-Invariant}
    by
    introducing 
    semi-equivariant transition cocycles,
    which simultaneously generalise 
    Wigner's corepresentations
    \cite[pp.~334-335]{Wigner.1959.Group-theory-and-its-applications-to-the-quantum-mechanics-of-atomic-spectra}
\cite[pp.~169-172]{Jansen-Boon.1967.Theory-of-finite-groups}
    and
    equivariant transition cocycles.
    Using the results of 
    \cite{Kitson.2020.A-Semi-equivariant-Dixmier-Douady-Invariant},
    the obstruction 
    to the existence 
    of an orientifold Dirac operator 
    can be identified as 
    a semi-equivariant Dixmier-Douady class.
    The main results are as follows
\begin{enumerate}
\item 
    Definition \ref{defn:spinkstructure}
    defines 
    \emph{$\Spink$-structures}.
    These are the appropriate 
    analogue of $\Spinc$-strucutre
    for orientifolds.
\item 
    Definition \ref{defn:thirdOSWC}
    defines the
    \emph{third integral 
    orientifold Stiefel-Whiney class}, 
    denoted $W_3^{(\Gam,\eps)}$.
\item
    Corollary \ref{cor:obs}
    shows that
    $W_3^{(\Gam,\eps)}$
    is the obstruction 
    to the existence of a
    $\Spink$-structure.
\item 
    Corollary \ref{cor:clas}
    shows that
    $\Spink$-structure are
    classified by
        $H_\Gam^2(X,(\Z,\iot_\eps))$,
    the elements of which 
    correspond to
    semi-equivariant principal 
    $(\Unitary(1),\kap_\eps)$-bundles.
\item
    Theorem \ref{thm:SOUcond}
    provides an alternative criteria
    for the existence of a $\Spink$-structure
    $P \rightarrow \Fr(V)$
    based on 
    the existence of a
    semi-equivariant 
    $(\Unitary(1),\kap_\eps)$-bundle
    that
    compliments
    the equivariant frame bundle $\Fr(V)$.
\item 
    Definition \ref{defn:ODO}
    defines the (Clifford-linear)
    \emph{orientifold Dirac operator}
    and
    \emph{reduced orientifold Dirac operator}.
\end{enumerate}

These results yeild 
    the primary theorem of this paper.
\begin{thm*}[\ref{thm:RealDirac}]
Let $X$ be an orientifold 
    with orientifold group $(\Gam,\eps)$.
\begin{enumerate}
\item 
    If $W^{(\Gam,\eps)}_3(X)=0$, then 
    $X$ carries 
    an orientifold Dirac operator.
\item 
    If $W^{(\Gam,\eps)}_3(X)=0$
    and $\dim(X) = 8$, 
    then 
    $X$ carries 
    a reduced orientifold Dirac operator.
\end{enumerate}
In particular,
    if $X$ is an 
    $8$-dimensional Real manifold
    and $W^{(\Z_2,\id)}_3(X)=0$, 
    then $X$ carries a 
    reduced Real Dirac operator.
\end{thm*}
These results appeared originally
    in the authors thesis
    \cite{Kitson.2020.Dirac-operators-on-orientifolds}.
The relationship between 
    the constructions described in this paper
    and other work in the literature
    will be discussed in \S\ref{sec:RWAPP},
    along with some potential applications.

\section{Orientifolds}
\label{ch:orientifolds}

This section begins 
    with a discussion of 
    orientifold groups. 
    Orientifold groups are 
    topological groups equipped
    with a small amount of extra structure
    that allows them to act
    in a linear/anti-linear manner.
    The representation theory of such actions
    on finite dimensional 
    complex vector spaces
    can be reduced to the theory of 
    unitary representations 
    that are invariant 
    under a conjugate structure
    on the space of equivalence classes 
    of representations.
    This reduction is achieved 
    by using the Wigner's notion 
    of a corepresentation
    \cite[pp.~334-335]{Wigner.1959.Group-theory-and-its-applications-to-the-quantum-mechanics-of-atomic-spectra}
\cite[pp.~169-172]{Jansen-Boon.1967.Theory-of-finite-groups},
    which coincides precisely
    with that of a semi-equivariant
    $(\Unitary(n),\kap_\eps)$-valued 
    transition cocycle 
    \cite[\S3]{Kitson.2020.A-Semi-equivariant-Dixmier-Douady-Invariant}
    over a point.

After briefly defining \emph{orientifolds},
    \emph{orientifold bundles} 
    will be introduced
    as complex vector bundles equipped
    with linear/anti-linear actions.
    On any orientifold bundle,
    it is possible to construct 
    a hermitian metric
    that is compatible
    with the linear/anti-linear action.
    Moreover, the frame bundle 
    of an orientifold bundle 
    is a 
    $\Gam$-semi-equivariant 
    principal $(\Unitary(n),\kap_\eps)$-bundle. 
    Thus, a neat generalisation is formed, 
    in which    
    an orientifold bundle over a point 
    is an orientifold representation, and 
    the semi-equivariant
    $(\Unitary(n),\kap_\eps)$-valued 
    transition cocycle 
    of its frame bundle
    is the corresponding
    corepresentation.
    From this perspective, 
    the semi-equivariant transition cocycles
    defined in 
    \cite[\S3]{Kitson.2020.A-Semi-equivariant-Dixmier-Douady-Invariant}
    are generalised corepresentations.
As with equivariant bundles,
    orientifold bundles admit 
    a number of natural operations.
    Semi-equivariant cocycles
    again prove useful,
    in \S\ref{sec:opsonOBs},
    for defining and working with
    these operations.

\subsection{Orientifold Groups}
\label{sec:OriGroups}

Any group $\Gam$ which acts 
    by a combination of linear and anti-linear
    operators must have an index-$2$ subgroup 
    of elements which act via linear operators,
    and a complementary subset of elements 
    which act via anti-linear operators.
    In general, if $\Gam$ contains more than one 
    subgroup of index $2$,
    then the set of 
    orientifold representations of $\Gam$ 
    depends on 
    which of these groups is chosen as $\Gam^+$.
    These facts motivate the definition 
    of an orientifold group.
\begin{defn}
An \emph{orientifold group} $(\Gam,\eps)$ 
    is a Lie group 
    equipped with a non-trivial homomorphism
        $\eps: \Gam \rightarrow \Z_2$.
    For any orientifold group define
        $\Gam^+ := \ker(\eps)$ and 
        $\Gam^- := \Gam \setminus \ker(\eps)$. 
\end{defn}
\begin{defn}
A \emph{homomorphism} 
        $\vphi: (\Gam',\eps') \rightarrow (\Gam,\eps)$
    of orientifold groups
    is a group homomorphism such that 
    $\eps \circ \vphi = \eps'$.
\end{defn}
The next lemma collects some
    basic facts about orientifold groups.
\begin{lem}
If $(\Gam,\eps)$ is an orientifold group, then
\begin{enumerate}
    \item $\Gam^+ \subset \Gam$ is a normal subgroup
    \item $\Gam/\Gam^+ \iso \Z_2$
    \item  
                $
                    1 
                        \rightarrow 
                    \Gam^+ 
                        \rightarrow 
                    \Gam 
                        \os{\eps}{\rightarrow} 
                    \Z_2 
                        \rightarrow 
                    1
                $
            is an extension of topological groups
    \item $\gam^2 \in \Gam^+$ for all $\gam \in \Gam$
    \item $\Gam = \Gam^+ \dunion \Gam^- = \Gam^+ \dunion \ze\Gam^+$ 
                for any 
              $\ze \in \Gam^-$.
\end{enumerate}
\end{lem}

The simplest non-trivial example 
    of an orientifold group
    is provided by 
        $\id: \Z_2 \rightarrow \Z_2$.
    Given an orientifold group, 
    its semi-direct product with a $\Gam$-group 
    can yield another orientifold group. 
\begin{lem}
\label{exam:OGprodG}
Let $\eps: \Gam \rightarrow \Z_2$ be an orientifold group
    and $(G,\tht)$ be a $\Gam$-group.
    Then the group extension
    \begin{equation*}
    \xymatrix@R.5em{
                1\ar[r] & 
                \Gam^+ \sdp_\tht G \ar[r]^(.5){i} & 
                \Gam \sdp_\tht G   \ar[r]^(.6){\eps \circ \pi_1} & 
                \Z_2\ar[r] & 
                1 
                \\
                & (\gam,g)\ar@{|->}[r] & (\gam,g)\ar@{|->}[r] & \eps(\gam) &
             }
    \end{equation*}
    makes $\Gam \sdp_\tht G$ into an orientifold group.
    The notation
        $(\Gam,\eps) \sdp_{\tht} G$
    will be used to denote orientifold groups of this form.
\end{lem}

The following example commonly arises 
    when $G$ is a group of linear operators 
    and $\kap$ represents conjugation with respect to a fixed basis.
\begin{exam}\label{exam:OGsdpZ2}
    Let $(G,\tht)$ be a $\Z_2$-group with unit $e$,
    then $(\Z_2,\id) \sdp_\tht G$ is an orientifold group
    \begin{equation*}
    \xymatrix@R.5em{
                1\ar[r] & 
                G\ar[r]^(.35)i & 
                \Z_2 \sdp_\tht G\ar[r]^(.6){\id \circ \pi_1} & 
                \Z_2\ar[r] & 
                1 \\
                & 
                g\ar@{|->}[r] & 
                (z,g)\ar@{|->}[r] & 
                z.&
             }
    \end{equation*}
Note that the element $(-1,e) \in \Gam^-$ is an involution,
    \begin{equation*}
        (-1,e)^2 = ((-1)^2,e (-1 e)) = (-1^2,e^2) = (+1,e).
    \end{equation*}
\end{exam}

It is also possible to construct 
    examples in which $\Gam^-$ 
    does not contain an involution.

\begin{exam}
The Weil group 
    \cite[\S XV]{Artin-Tate.2009.Class-field-theory}
    of $\R$
    is the subgroup
    $\C^\times \dunion \C^\times j 
        \subset \H^\times$
    of the multiplicative group 
    of quaternions.
    It
    fits
    into the non-split extension
    \begin{equation*}
    \xymatrix@R.5em{
                1
                    \ar[r] & 
                \C^\times 
                    \ar[r]^(.375){} & 
                \C^\times \dunion \C^\times j
                    \ar[r]^(.5){} & 
                \Gal(\C / \R)
                    \ar[r] & 
                1 
                \\
                & 
                & 
                j\ar@{|->}[r] & 
                -1&
             }
    \end{equation*}
    of $\C^\times$ 
    by $\Gal(\C / \R) \iso \Z_2$,
    making it into an orientifold group.
    Note that there is no element 
    $\ze \in \C^\times j = \Gam^-$
    such that $\ze^2 = 1$.
\end{exam}

\begin{exam}
If $H := \set{\pm 1, \pm i}$ is the orientifold group
    equipped with the homomorphism $q(h) := h^2$,
    then 
        $\Gam := (H,q) \sdp_\tht G$ 
    is the orientifold group
    \begin{equation*}
    \xymatrix@R.5em{
                1\ar[r] & 
                \set{\pm 1} \sdp_\tht G \ar[r]^(.45){i} & 
                \set{\pm 1,\pm i} \sdp_\tht G\ar[r]^(.65){q \circ \pi_1} & 
                \Z_2\ar[r] & 
                1 
                \\
                & (h,g)\ar@{|->}[r] & (h,g)\ar@{|->}[r] & h^2. &
             }
    \end{equation*}
    If $(h,g) \in \Gam^-$,
    then $h = \pm i$ and
        $(h,g)^2 = (h^2,g(hg)) = (-1,g(hg)) \in \Gam^+$.
    Thus, there is no element $\gam \in \Gam^-$ 
    such that $\gam^2 = (1,e)$.
\end{exam}

Given an orientifold group $(\Gam,\eps)$, 
    the parity information provided by $\eps$
    can be used when defining actions on various objects.
    Three different types 
    of actions of an orientifold group 
    will be distinguished. 
The first type of action uses the parity information 
    assigned to group elements 
    to dictate whether an element 
    acts linearly or anti-linearly.
It will be neccesary to define
    these actions on a variety
    of $\C$-modules 
    from different categories, 
    including
        complex vector spaces, 
        complex vector bundles, and 
        algebras over $\C$.
Given objects $X$ and $Y$ in an appropriate category, define
    \begin{align*}
        \Hom^{+1}(X,Y) &:= \Hom(X,Y)
        \\
        \Hom^{-1}(X,Y) 
            &:= 
        \set{ 
             a_{\conj{Y}} \circ \vphi 
                 \mid 
             \vphi \in \Hom(X,\conj{Y}) 
            },
    \end{align*}
    where 
        $a_{\conj{Y}}: \conj{Y} \rightarrow Y$
    is the identity map on the underlying set for $Y$.
    The map $a_{\conj{Y}}$ is anti-linear
    and the elements of 
        $\Hom^{-1}(X,Y)$ 
    can be considered as anti-linear homomorphisms.
    The conjugation map 
        $Y \mapsto \conj{Y}$ 
    changes the $\C$-module structure of $Y$
    to its conjugate $\C$-module structure, and,
    depending on the category, 
    it may change other structures on $Y$.
    For example, 
    the conjugate of a Hilbert space 
    carries a conjugate inner product.
    Denote the disjoint union 
    of $\Hom^+$ and $\Hom^-$ by $\Hom^\pm$.
    The spaces $\End^\pm$ and $\Aut^\pm$
    are defined similarly.

\begin{defn}
\label{defn:orientfoldAction}
Let $(\Gam,\eps)$ be an orientifold group.
    An \emph{orientifold action} is a homomorphism
        $\rho: \Gam \rightarrow \Aut^{\pm}(W)$
    such that
    \begin{equation*}
        \rho(\gam) \in \Aut^{\eps(\gam)}(W).
    \end{equation*}
\end{defn}

A second type of action
    uses an involution $\rho$ 
    to define an action of $\Gam$.
    Typically, 
    an involution of this type
    represents 
    the change of some structure 
    to a conjugate structure, 
    occuring in parallel 
    with the application of an orientifold action.

\begin{defn}
    An \emph{involutive action} 
    of an orientifold group,
    is an action of the form
    \begin{equation}
        \rho \circ \eps: \Gam \rightarrow \Z_2 \rightarrow \Aut(Y),
    \end{equation} 
    where
      $\rho: \Z_2 \rightarrow \Aut(Y)$
    is an involution.
\end{defn}

\begin{exam}
\label{exam:orientifoldActions}
Some examples of involutive actions are
\begin{enumerate}
    \item 
        $\iot^{p,q}_\eps: \R^{p,q} \rightarrow \R^{p,q}$,
        where $\R^{p,q} := \R^p \dsum \R^q$,
        $\iot^{p,q}: (x,y) \mapsto (x,-y)$.
    \item 
        $\kap_\eps: \GL(n,\C) \rightarrow \GL(n,\C)$,
        where $\kap_\eps$ is elementwise conjugation
        on the standard matrix representation of $\GL(n,\C)$.
    \item 
        $d\tht_\eps: \aL{g} \rightarrow \aL{g}$, 
        where $\aL{g}$ is a Lie algebra and
            $\tht: G \rightarrow G$
        is an involution on its Lie group.
\end{enumerate}
\end{exam}
Throughout this paper, it will be assumed
    that $\C$ is equipped with the orientifold
    action $\kap_\eps$.

Of course, the parity of the group elements can also be ignored.
    This type of action occurs on an orientifold and 
    its tangent bundle.
    In order to differentiate it from the other types of action, 
    it will be refered to as a \emph{basic} action.

\subsection{Orientifolds}
\label{sec:orientifolds}

In order to maintain a clear
    focus on anti-linear symmetry, 
    only the simplest definition 
    of an orientifold will be treated.
These orientifolds
    are essentially global quotient orbifolds 
    with a small amount of extra structure.
Using the language of \S\ref{sec:OriGroups},
    they could be described as
    manifolds equipped 
    with a basic action 
    of an orientifold group.
The origin 
    of the term orientifold 
    is in string theory,
    where
    orientifolds are often considered to have 
    a sign choices $\pm 1$ 
    associated to the connected components 
    of their fixed point sets.
    However, these sign choice structures will not be considered here.
\begin{defn}\label{defn:orientifold}
    An \emph{orientifold} 
    is a compact manifold $X$ 
    equipped with a smooth action 
    \begin{equation*}
        \rho: \Gam \rightarrow \Diff(X),
    \end{equation*}
    where $\Gam$ is 
    a finite orientifold group.
The category of orientifolds 
    with orientifold group 
        $\eps: \Gam \rightarrow \Z_2$ 
    will be denoted $\Ori_{(\Gam,\eps)}$.
\end{defn}

\begin{exam}
    Let $\Gam$ be any orientifold group.
    Then $\R^{p,q} := \R^p \dsum \R^q$ 
    equipped with the involutive action 
    induced by
            $(x,y) \mapsto (x,-y)$
        is an orientifold.
    This orientifold will be used 
    to form suspensions 
    in orientifold $K$-theory.
\end{exam}
\begin{exam}
    Let $X \in \Ori_{(\Gam,\eps)}$ with $\Gam$-action $\sig$.
    The tangent bundle $TX$ equipped 
    with the basic $\Gam$-action $d\sig$
    is again an orientifold.
    The $K$-theory of this orientifold is 
    the target space 
    of the $8$-fold Thom isomorphism 
    in orientifold $K$-theory. 
\end{exam}
The category of real vector bundles 
    equipped with a basic action of the orientifold group $(\Gam,\eps)$ 
    will be denoted $\Vect_{(\Gam,\eps)}(X,\R)$.
    The isomorphism classes of such bundles 
        will be denoted $\Vect_{(\Gam,\eps)}^\iso(X,\R)$.

\subsection{Orientifold Bundles}
\label{sec:OBundles}

Orientifold bundles 
    are the main object of interest 
    in the study of orientifolds.
    In the language of \S\ref{sec:OriGroups},
    they are complex vector bundles 
    carrying orientifold actions
    that cover the action 
    on the base orientifold.

\begin{defn}
If $\pi: E \rightarrow X$ is a complex vector bundle,
    define
        $\Aut_{\Diff}(E)$
    to be the set of maps
        $\vphi: E \rightarrow E$ 
    such that 
    \begin{enumerate}
    \item 
        $\pi \circ \vphi(e) = f \circ \pi(e)$,
        for some $f \in \Diff(X)$
        and all $e \in E$.
    \item
        $\vphi: E_x \rightarrow E_{f(x)}$
        is a linear bijection, for all $x$.
    \end{enumerate}
\end{defn}
    
\begin{defn}
An \emph{orientifold bundle} $\pi: E \rightarrow X$ 
    is a complex vector bundle
    equipped with an orientifold action
    \begin{equation*}
        \tau: \Gam \rightarrow \Aut^\pm_{\Diff}(E)
    \end{equation*}
    such that $\pi(\gam v) = \gam \pi(v)$. 
\end{defn}

The category of 
    orientifold bundles 
    over $X \in \Ori_{(\Gam,\eps)}$ 
    will be denoted $\Vect_{(\Gam,\eps)}(X,\C)$.
    The set of isomorphism classes of orientifold bundles 
    will be denoted $\Vect_{(\Gam,\eps)}^\iso(X,\C)$.
\begin{exam}
A linear/anti-linear representation $(V,\rho)$
    is an orientifold bundle over a point.
    Such a bundle will also be refered to as 
    an \emph{orientifold representation}.
    If $(X,\sig)$ is an orientifold 
    and $(V,\rho)$ 
    is an orientifold representation,
    then an orientifold bundle of the form
    \begin{equation*}
        (X \times V,\sig \times \rho)
    \end{equation*}
    will be described as 
    a \emph{trivial} orientifold bundle.
Note that if $\eps$ is non-trivial, 
    then every orientifold bundle 
    for $(\Gam,\eps)$
    carries at least one anti-linear map, 
    and so there is 
    no orientifold bundle $(E,\tau)$ 
    such that 
    $\tau_\gam = \id$ for all $\gam \in \Gam$. 
\end{exam}

Just as in the equivariant setting, 
    it is possible 
    to average an hermetian metric 
    on an orientifold bundle
    to make it compatible 
    with the orientifold action.
    The averaging process needs 
    to be twisted with conjugation 
    to account 
    for the anti-linearity of the action,
    as does the compatibility condition.

\begin{defn}
    An \emph{orientifold metric} 
    on an orientifold bundle $E$
    is an hermitian metric $h$ on $E$ such that, 
    for all $v_1,v_2 \in E$ and $\gam \in \Gam$,
    \begin{equation*}
        h(\gam v_1, \gam v_2)_{\gam x} = \gam h(v_1,v_2)_x.
    \end{equation*}
\end{defn}

\begin{prop}\label{prop:IHM}
Every orientifold vector bundle $E$
    over a paracompact orientifold $X$ 
    carries an orientifold metric.
\begin{proof}
It is a standard result that 
    every complex vector bundle 
    over a paracompact space 
    carries an hermitian metric
        \cite[Lemma 2]{Swan.1962.Vector-bundles-and-projective-modules}.
Given an hermitian metric $h$
        on an orientifold bundle $E$,
        define
        \begin{equation*}
            h_\Gam(u,v)_x 
              = 
            \sum_{\gam \in \Gam} \gam^\inv h(\gam u,\gam v)_{\gam x}.
        \end{equation*}
    This metric is an orientifold metric as
    \begin{align*}
        h_\Gam(\gam u, \gam v)_{\gam x} 
          &= 
        \sum_{\gam' \in \Gam} 
            \gam'^\inv h(\gam' \gam u,\gam' \gam v)_{\gam' \gam x} 
        \\&\qquad
          = 
        \sum_{\gam'' := \gam'\gam \in \Gam} 
            \gam\gam''^\inv h(\gam'' u, \gam'' v)_{\gam'' x} 
        \\&\qquad\qquad
          = 
        \gam \sum_{\gam'' \in \Gam} 
            \gam''^\inv h(\gam'' u, \gam'' v)_{\gam'' x} 
          = 
        \gam h_\Gam(u,v)_x.
    \end{align*}
\end{proof}
\end{prop}

Using an orientifold metric, 
    it is possible to split
    sequences of orientifold bundles.

\begin{cor}
    Let $X$ be a paracompact orientifold.
    If 
    \begin{equation*}
        0
          \rightarrow 
        E' 
            \os{\vphi'}{\rightarrow} 
        E 
            \os{\vphi}{\rightarrow} 
        E''
    \end{equation*}
    is an exact sequence of orientifold bundles over $X$,
    then 
        $E \iso E' \dsum E''$.
\begin{proof}
    By Proposition \ref{prop:IHM}, there exists an orientifold metric $h$ on $E$.
    It is a standard result 
        \cite[Proposition 2]{Swan.1962.Vector-bundles-and-projective-modules}
    that $h$ determines a projection $p: E \rightarrow E$
    and a splitting of complex vector bundles
        $E = \im(p) \dsum \ker(p) \iso E' \dsum E''$.
    The projection $p$ is defined fibrewise by 
    \begin{align*}
        p_x: E_x &\rightarrow E_x \\
         v &\mapsto \sum_i \frac{ h(v,b_i)_x }{ h(b_i,b_i)_x } b_i,
    \end{align*}
    where $\set{b_i}$ is any basis for $\vphi'(E')_x$.
    Therefore, if $p_x(v)=0$, 
    then $h(v,b_i)_x = 0$ for all $i$, 
    and
    \begin{align*}
         p_{\gam x}(\gam v) 
            &= 
        \sum_i 
            \frac{h(\gam v,\gam b_i)_{\gam x}}
                 {h(\gam b_i,\gam b_i)_{\gam x}} 
                 (\gam b_i)
            \\&\qquad= 
        \sum_i 
            \frac{\gam h(v,b_i)_{x}}
                 {\gam h(b_i,b_i)_{x}} 
                 (\gam b_i)
            = 
        \sum_i 
            \frac{\gam 0}
                 {\gam h(b_i,b_i)_{x}} 
                 (\gam b_i)
            = 0.
    \end{align*}
    Thus, 
    $\ker(p)$ is invariant 
    under the action of $\Gam$, 
    as is the given splitting.
\end{proof}
\end{cor}


Next, the frame bundle 
    of an orientifold bundle will be examined.

\begin{defn}
The \emph{frame bundle} $\Fr(E)$ 
    of an orientifold bundle $E$ 
    is the principal $\GL(n,\C)$-bundle of frames
    for the total space of $E$,
    equipped with a left $\Gam$-action defined on a frame 
    $s = (s_1, \ldots, s_n) \in \Fr(E)_x$ 
    by
        $(\gam s)_i = \gam s_i$.
\end{defn}

Although the frame bundle of an orientifold 
    is defined in the same manner
    as that of an equivariant bundle, 
    the anti-linearity present 
    in the $\Gam$-action gives it different properties.
    In particular, there is a mild noncommutivity 
    between 
      the left action of $\Gam$ and 
      the right action of the structure group $\GL(n,\C)$.
    This non-commutivity makes 
    the frame bundle of an orientifold bundle 
    into 
    a $\Gam$-semi-equivariant 
    principal $(\GL(n,\C),\kap_\eps)$-bundle
    \cite[\S2]{Kitson.2020.A-Semi-equivariant-Dixmier-Douady-Invariant}.
\begin{prop}
\label{prop:OBFB}
Let $E$ be an orientifold bundle and 
    consider $\GL(n,\C)$ to be equipped 
    with the involutive action 
    of $(\Gam,\eps)$ induced by conjugation.
    Then, 
    \begin{equation*}
        \Fr(E;\GL(n,\C)) \in \PBundles_{\Gam}(X,(\GL(n,\C),\kap_\eps)).
    \end{equation*}
    In particular,
    the left and right actions 
    on the frame bundle 
    satisfy 
    \begin{equation*}
      \gam(sg) = (\gam s)(\gam g),
    \end{equation*}
    for 
        $\gam \in \Gam$, 
        $s \in \Fr(E)$ and 
        $g \in \GL(n,\C)$.
  \begin{proof}
    The action of $g$ on a frame $s$ is given by
      $(sg)_j = \sum_{1 \leq i \leq n} s_ig_{ij}$.
    Thus,
    \begin{align*}
        \gam(sg)_j 
          &= \sum_{1 \leq i \leq n} \gam(s_ig_{ij}) 
          = \sum_{1 \leq i \leq n} (\gam s_i)(\gam g_{ij}) 
          = \sum_{1 \leq i \leq n} (\gam s)_i(\gam g)_{ij} 
    \\&\qquad
          = ((\gam s)(\gam g))_j.
    \end{align*}
  \end{proof}
\end{prop}

Note that, by using an orientifold metric,
    the structure group can always 
    be reduced to $(\Unitary(n),\kap_\eps)$,
    where $\kap_\eps$ is the action induced on $\Unitary(n)$
    by its inclusion into $\GL(n,\C)$.

\subsection{Operations on Orientifold Bundles}
\label{sec:opsonOBs}

Some basic operations on orientifold bundles
    will now be defined.
    It will be useful to make these definitions 
    in terms of semi-equivariant cocycles
    \cite[\S 3]{Kitson.2020.A-Semi-equivariant-Dixmier-Douady-Invariant}. 
    To start with, 
    consider the following operations on $\Gam$-groups.

\begin{defn}\label{def:GLn}
Let $a^k \in \GL(\C^{m_k})$,
    and denote by $[a_{ij}]$ the matrix representation of an element $a \in \GL(\C^m)$
    with respect to the standard basis of $\C^m$.
    Define the following operations
    \begin{enumerate}
        \item 
                The \emph{dual} 
                        $a^* \in \GL(\C^m)$,
                    \begin{equation*}
                        [(a^*)_{ij}]
                            :=
                        ([a_{ij}]^t)^\inv
                    \end{equation*}
        \item 
                The \emph{direct sum} 
                        $a^1 \dsum a^2 \in \GL(\C^{m_1 + m_2})$,
                    \begin{equation*}
                        [(a^1 \dsum a^2)_{ij}]
                            :=
                        \begin{pmatrix}
                         [a^1_{ij}]  &   0              \\
                        0          &   [a^2_{ij}] 
                        \end{pmatrix}.
                    \end{equation*}
        \item 
                The \emph{tensor product}
                            $a^1 \ten a^2 \in \GL(\C^{m_1m_2})$,
                        \begin{equation*}
                            [(a^1 \ten a^2)_{ij}]
                                := 
                            \begin{pmatrix}
                             a^1_{11}[a^2_{ij}]    & \ldots & a^1_{1m}[a^2_{ij}]    \\
                            \vdots & \ddots & \vdots \\
                             a^1_{m1}[a^2_{ij}]    & \ldots & a^1_{mm}[a^2_{ij}] 
                            \end{pmatrix}.
                        \end{equation*}
    \end{enumerate}
\end{defn}

Examining Definition \ref{def:GLn}, 
    it is clear that the dual, direct sum and tensor product on the groups $\GL(\C^m)$ 
    are compatible with involutive $\Gam$-actions induced by conjugation.

\begin{lem}\label{lem:tenGLn}
The dual, direct sum and tensor product operations are homomorphisms
    \begin{align*}
        *
            &: 
        (\GL(\C^m),\kap_\eps) 
            \rightarrow 
        (\GL(\C^m),\kap_\eps) 
        \\
        \dsum
            &: 
        (\GL(\C^{m_1}),\kap_\eps) \times (\GL(\C^{m_2}),\kap_\eps) 
            \rightarrow 
        (\GL(\C^{m_1 + m_2}),\kap_\eps) 
        \\
        \ten
            &: 
        (\GL(\C^{m_1}),\kap_\eps) \times (\GL(\C^{m_2}),\kap_\eps) 
            \rightarrow 
        (\GL(\C^{m_1m_2}),\kap_\eps)
    \end{align*}
    of $\Gam$-groups.
\end{lem}     
    
Lemma \ref{lem:tenGLn}, allows the dual, direct sum and tensor product 
    of $(\GL(m,\C),\kap_\eps)$-valued transition cocycles 
    \cite[\S3]{Kitson.2020.A-Semi-equivariant-Dixmier-Douady-Invariant}
    to be defined in the obvious way.
    Pullbacks of cocycles can also be defined.
    It is routine to prove that these satisfy 
    the semi-equivariant cocycle condition.

\begin{defn}
Let $\phi^i \in 
        \TCocycles_{(\Gam,\eps)}
            (\cU,X,(\GL(\C^{m_i}),\kap_\eps))$
    and $f: X \rightarrow Y$ be a homomorphism orientifolds. 
    The 
        \emph{pullback},
        \emph{dual},
        \emph{direct sum}, and 
        \emph{tensor product}
    are defined, respectively,
    by
    \begin{alignat*}{3}
    (f^*\phi)_{ba}(\gam, x) 
    &:= 
    \phi_{ba}(\gam, f(x)) 
    &&
    \in \TCocycles_{(\Gam,\eps)}(f^*\cU, Y,(\GL(\C^m),\kap_\eps))
    \\
    (\phi^*)_{ba}(x,\gam) 
    &:= 
    \phi_{ba}(x,\gam)^*
    &&
    \in \TCocycles_{(\Gam,\eps)}(\cU,X,(\GL(\C^m),\kap_\eps))
    \\
    (\phi^1 \dsum \phi^2)_{ba}(x,\gam) 
    &:= 
    \phi^1_{ba}(x,\gam) \dsum \phi^2_{ba}(x,\gam)
    &&
    \in \TCocycles_{(\Gam,\eps)}(\cU,X,(\GL(\C^{m_1+m_2}),\kap_\eps))
    \\
    (\phi^1 \ten \phi^2)_{ba}(x,\gam) 
    &:= 
    \phi^1_{ba}(x,\gam) \ten \phi^2_{ba}(x,\gam)
    &&
    \in
    \TCocycles_{(\Gam,\eps)}(\cU,X,(\GL(\C^{m_1m_2}),\kap_\eps)),
    \end{alignat*}
    where 
    $f^*\cU := \set{ f^\inv(U_a) \mid a \in A }$     is the pullback of 
    $\cU := \set{U_a \mid a \in A}$.
\end{defn}

The above operations on cocycles 
    induce operations on orientifold bundles 
    via the 
    semi-equivariant 
    associated bundle construction, 
    see Definition \ref{defn:SEAB}.

\begin{defn}
Let $E_i \in \Vect_{(\Gam,\eps)}^{m_i}(X,\C)$.
    Let $\phi^i$ 
    denote a semi-equivariant cocycle 
    associated $\Fr(E_i)$ 
    by 
    \cite[Prop.~12]{Kitson.2020.A-Semi-equivariant-Dixmier-Douady-Invariant}, and 
    $P^\phi$ denote 
    the semi-equivariant 
    principal bundle constructed 
    from a cocycle $\phi$ 
    via 
    \cite[Prop.~15]{Kitson.2020.A-Semi-equivariant-Dixmier-Douady-Invariant}.
    $P^\phi$ denote 
    The
        \emph{pullback},
        \emph{dual},
        \emph{direct sum}, and 
        \emph{tensor product}
    operations on orientifold bundles
    are defined, respectively, by
    \begin{alignat*}{3}
    f^*E 
    &:= 
    P^{f^*\phi} \times_{(\GL(m,\C),\kap_\eps)} (\C^{m},\kap_\eps) 
    &&
    \in \Vect_{(\Gam,\eps)}^m(X,\C)
    \\
    E^* &:= P^{\phi^*} \times_{(\GL(m,\C),\kap_\eps)} ((\C^{m})^*,\kap_\eps)
    &&
    \in
    \Vect_{(\Gam,\eps)}^m(X,\C)
    \\
    E_1 \dsum E_2 
    &:= 
    P^{\phi_1 \dsum \phi_2} \times_{(\GL(m_1+m_2,\C),\kap_\eps)} (\C^{m_1+m_2},\kap_\eps)
    &&
    \in \Vect_{(\Gam,\eps)}^{m_1+m_2}(X,\C)
\\
    E_1 \ten E_2 
    &:= 
    P^{\phi_1 \ten \phi_2} \times_{(\GL(m_1m_2,\C),\kap_\eps)} (\C^{m_1m_2},\kap_\eps)
    &&
    \in \Vect_{(\Gam,\eps)}^{m_1m_2}(X,\C),
    \end{alignat*}
    where $\kap_\eps:(\C^m)^* \rightarrow (\C^m)^*$ is the action defined by
        $(\gam \lam)(z) := \gam \lam(\gam^\inv z)$.
\end{defn}

As in the non-equivariant setting, 
    it is possible to construct the bundle
    of homomorphisms between 
    two orientifold bundles 
    using their tensor products and duals.

\begin{prop}
Let $E_i \in \Vect_{(\Gam,\eps)}^{m_i}(X,\C)$.
    Homomorphisms
    in
    $\Hom(E_1,E_2)$ 
    correspond bijectively
    to 
    equivariant sections
    of the orientifold bundle 
    $E_2 \ten E_1^*$.
\end{prop}


\subsection{Classification of $\Spinc$-Structures on Orientifolds}

\label{sec:OSS}

In order to define and classify $\Spinc$-structures for orientifolds,
    it is neccesary to consider 
    the interaction of Clifford algebras and 
    the $\Spin$ groups with orientifold actions.
    The idea is to complexify results which apply to real Clifford algebras, 
    whilst keeping track of the associated conjugation maps.
    These maps can then be used 
    to define involutive actions of orientifold groups.
To begin, the definitions of 
    the real Clifford algebra, $\Spin$ group, and adjoint map are recalled.

\begin{defn}
    The Clifford algebra $\Cl_n$ is the algebra generated
    by the standard basis $\set{e_i}$ of $\R^n$
    subject to the relations
    $e_i^2 = -1$ and $e_ie_j + e_je_i = 0$.
\end{defn}
Note that the set
        $\set{e_{i_1} \cdots e_{i_k} \in \Cl_n \mid i_1 < \cdots < i_k}$
    is a basis for $\Cl_n$.
The group $\Spin(n)$ sits inside $\Cl_n$.
    Elements of $\Spin(n)$ are products of an even number of unit vectors from $\R^n$.

\begin{defn}
    The group $\Spin(n)$ is defined by 
    \begin{equation*}
        \Spin(n) 
            := 
        \set{ x_1 \cdots x_{2k} \mid x_i \in \R^n, \Norm{x_i} = 1 } 
            \subset \Cl_n.
    \end{equation*}
\end{defn}

If $g \in \Spin(n)$ and $x \in \R^n$ one can show that $g x g^\inv \in \R^n$.
    The transformation 
        $x \mapsto gxg^\inv$
    defines an element of $\SO(n)$, and 
    the resulting assignment
        $\Spin(n) \rightarrow \SO(n)$
    is a double covering.

\begin{defn}
    The \emph{adjoint map} 
        $\Ad: \Spin(n) \rightarrow \SO(n)$ 
    is defined, for $g \in \Spin(n)$, $x \in \R^n$, by 
    \begin{equation*}
        \Ad_g(x) := gxg^\inv.
    \end{equation*}
\end{defn}

For applications to orientifolds, 
    it is neccesary to work with the
    complexifications of $\Cl_n$ and $\Spin(n)$.
    These complexifications are equipped with conjugation maps
    which induce involutive actions of orientifold groups.
    The complexified adjoint map is a homomorphism of $\Gam$-groups.
\begin{defn}
\label{defn:ClcEtc}
Let $(\Gam,\eps)$ be an orientifold group and 
    define the following
    \begin{enumerate}
    \item 
        $(\Clc_n,\kap_\eps) := \Cl_n \ten \C$
        with the $\Gam$-action 
        $
          \kap_\eps(\vphi \ten z) 
              := 
          \vphi \ten \kap_\eps(z)
        $
    \item 
        $
          (\Spinc(n),\kap_\eps) 
              := 
          ( \Spin(n) \times \Unitary(1) ) 
              / 
          \set{\pm(1,1)}
        $
        with the induced action 
        $\kap_\eps[g,z] := [g,\kap_\eps(z)]$
    \item 
        $
         \Ad^c: (\Spinc(n),\kap_\eps) 
            \rightarrow 
         (\SO(n),\id_\eps)$ 
        defined by 
        $\Ad^c[g,z] := \Ad(g)$.
    \end{enumerate}
\end{defn}
Note that $\Ad^c \circ \kap_\eps[g,z] = \Ad^c[g,z]$.
The properties of $\Ad^c$, and 
    the decomposition of $\Spinc(n)$, 
    produce 
    two central exact sequences 
    of $\Gam$-groups 
    about $\Spinc(n)$.
    These sequences fit into 
    the following diagram
    \begin{equation}
    \label{fig:SpincEx}
    \begin{aligned}
    \end{aligned}
    \begin{aligned}
    \tikz[xscale=3,yscale=1.5]{\tikzstyle{every node}=[font=\small,scale=1]
       \node(A) at (1,2)    {1};
       \node(B) at (2,2)    {$(\Spin(n),\id_\eps)$}
            edge [<-] node[auto] {} (A);
       \node(C) at (3,2)    {$(\Spinc(n),\kap_\eps)$}
            edge [<-] node[auto] {} (B);
       \node(D) at (4,2)    {$(\Unitary(1),\kap_\eps)$}
            edge [<-] node[auto] {} (C);
       \node(E) at (4.5,2)    {1}
            edge [<-] node[auto] {} (D);

       \node(F) at (3,0)    {1};
       \node(G) at (3,1)    {$(\SOrth(n),\id_\eps)$}
            edge [->] node[auto] {} (F)
            edge [<-] node[auto,swap] {$\Ad^c$} (C)
            edge [<-] node[auto] {$\Ad$} (B);
       \node(H) at (3,3)    {$(\Unitary(1),\kap_\eps)$}
            edge [->] node[auto] {} (C)
            edge [->] node[auto,font=\small] {$q$} (D);
       \node(I) at (3,4)    {1}
            edge [->] node[auto] {} (H);

       \node(J) at (0.5,3.5)    {1};
       \node(K) at (1,3)    {$(\Z_2,\id_\eps)$}
            edge [<-] node[auto] {} (J)
            edge [->] node[auto] {} (B);
       \node(L) at (4,0)    {1}
            edge [<-] node[auto] {} (G);

       \node(M) at (1,5)    {1};
       \node(N) at (2,4)    {$(\Z_2,\id_\eps)$}
            edge [->] node[auto] {} (H)
            edge [<-] node[auto] {} (M);
       \node(O) at (4.5,1.5)    {1}
            edge [<-] node[auto] {} (D);
    },
    \end{aligned}
    \end{equation}
    where $q$ is the square map.
    The above sequences will be used 
    to classify $\Spinc$-structures for orientifolds.

Having examined 
    semi-equivariance, 
    orientifolds, and 
    orientifold actions on $\Spinc(n)$,
    it is now possible to define a notion of $\Spinc$-structure which is
    appropriate for orientifolds.
\begin{defn}
\label{defn:spinkstructure}
    An \emph{$\Spink$-structure} 
    for a real $\Gam$-equivariant 
    vector bundle $V$ 
    over an orientifold
    is a semi-equivariant lifting 
    $\vphi: P \rightarrow \Fr(V)$ 
    by
    $\Ad^c: (\Spinc(n),\kap_\eps) \rightarrow (\SOrth(n),\id_\eps)$.
\end{defn}
If $V$ has a $\Spink$-structure, 
    then it is said to be 
    \emph{$\Spink$-oriented}.
    If the tangent bundle $TM$ 
    of an orientifold $M$ is $\Spink$-oriented,
    then $M$ is said to be 
    $\Spink$-\emph{oriented}.

The $\Spink$-structures associated 
    to a vector bundle $V$
    can be classified using the results of 
    \cite{Kitson.2020.A-Semi-equivariant-Dixmier-Douady-Invariant}.
    The following theorem is obtained by 
    applying 
    \cite[Theorem 41]{Kitson.2020.A-Semi-equivariant-Dixmier-Douady-Invariant}
    to the central exact sequence running vertically 
    in diagram \eqref{fig:SpincEx}.

\begin{thm}\label{thm:OSpincSeq}
The central exact sequence
    \begin{equation*}
        1 
            \rightarrow  
        (\Unitary(1),\kap_\eps) 
            \rightarrow 
        (\Spinc(n),\kap_\eps) 
            \os{\Ad^c}{\rightarrow} 
        (\SOrth(n),\id_\eps)
            \rightarrow 
        1,
    \end{equation*}
    induces an exact sequence
    \begin{align*}
        H^1_\Gam(X,(\Unitary(1),\kap_\eps))     
          &\os{}{\rightarrow}
        \TCocycles^1_\Gam(X,(\Spinc(n),\kap_\eps))
          \os{\Ad^c}{\rightarrow}
          \ldots 
    \\\qquad
        &\qquad\ldots 
        \os{\Ad^c}{\rightarrow}
        \TCocycles^1_\Gam(X,(\SOrth(n),\id_\eps))
          \os{\Del_{sc}}{\rightarrow}
        H^2_\Gam(X,(\Unitary(1),\kap_\eps)).
    \end{align*}
\end{thm}
Theorem \ref{thm:OSpincSeq} has the following corollaries,
    which classify $\Spink$-structures 
    in terms of semi-equivariant cohomology 
    with coefficients in $(\Unitary(1),\kap_\eps)$.
\begin{cor}
\label{cor:obs}
A real $\Gam$-equivariant vector bundle $V$ 
    over an orientifold
    has a $\Spink$-structure 
    if and only if 
    $\Del_{sc}(\phi^V) = 1$,
    where $\phi^V$ is the transition cocycle for $V$.
\end{cor}

\begin{cor}
\label{cor:clas}
A given $\Spink$-structure is unique 
    up to tensoring by 
    semi-equivariant principal 
    $(\Unitary(1),\kap_\eps)$-bundles.
\end{cor}

To obtain an obstruction class with integer coefficients, 
    involutive actions can be taken on the groups 
    in the exponential exact sequence.
    This results in the following proposition.
\begin{lem}
\label{lem:expExCo}
The exponential exact sequence 
    \begin{equation}\label{eq:expEx}
        0 \rightarrow   
        (\Z,\iot_\eps) \rightarrow
        (\R,\iot_\eps) \os{\exp}{\rightarrow}
        (\Unitary(1),\kap_\eps) \rightarrow
        1
    \end{equation}
    induces isomorphisms
    \begin{equation*}
        H^p_\Gam(X,(\Unitary(1),\kap_\eps))
          \os{\Del_{\exp}^p}{\rightarrow}
        H^{p+1}_\Gam(X,(\Z,\iot_\eps)),
    \end{equation*}
    where $\iot_\eps$
    is the involutive orientifold action 
    induced by the map $t \mapsto -t \in \R$.
\begin{proof}
By 
    \cite[Theorem 38]{Kitson.2020.A-Semi-equivariant-Dixmier-Douady-Invariant},
    the exact sequence \eqref{eq:expEx} induces a long exact sequence
    \begin{equation*}
        H^p_\Gam(X,(\Z,\iot_\eps))     
          \os{}{\rightarrow}
        H^p_\Gam(X,(\R,\iot_\eps))
          \os{\exp}{\rightarrow}
        H^p_\Gam(X,(\Unitary(1),\kap_\eps))
          \os{\Del_{\exp}^p}{\rightarrow}
        H^{p+1}_\Gam(X,(\Z,\iot_\eps)).
    \end{equation*}
    The cohomology groups 
    $H^p_\Gam(X,(\R,\iot_\eps))$ vanish 
    for all $p$,
    due to the existence of 
    a smooth partition of unity on $X$.
    Therefore, 
    the maps $\Del_{\exp}^p$ are isomorphisms.
\end{proof}
\end{lem}

Using Lemma \ref{lem:expExCo}, it is possible to define an analogue
    of the third integral Stiefel-Whiney class.

\begin{defn}
\label{defn:thirdOSWC}
The \emph{third integral orientifold Stiefel-Whiney class}
    is defined by
    \begin{equation*}
        W_3^{(\Gam,\eps)}(V) 
            := 
        \Del_{exp} \circ \Del_{sc}(\phi^V) 
            \in 
        H_\Gam^3(X,(\Z,\iot_\eps)),
    \end{equation*}
    where $\phi^V$ is the transition cocycle associated to $V$.
\end{defn}

Corollaries \ref{cor:obs} and \ref{cor:clas} 
    can then be restated
    in terms of semi-equivariant cohomology 
    with coefficients in $(\Z,\iot_\eps)$.

\begin{cor}
\label{cor:W3}
A real $\Gam$-equivariant bundle $V$ is $\Spink$-oriented
    if and only if $W_3^{(\Gam,\eps)}(V) = 0$. 
\end{cor}

\begin{cor}
The $\Spink$-structures on 
    a $\Spink$-oriented 
    real $\Gam$-equivariant vector bundle 
    are in bijective correspondence 
    with the elements of 
        $H_\Gam^2(X,(\Z,\iot_\eps))$.
\end{cor}

It is possible to further isolate 
    the semi-equivariance in a $\Spink$-structure by 
    splitting it via the decomposition
    \begin{equation*}
        (\Spinc(n),\kap) 
            \iso 
        (\SOrth(n),\id) \times_{\Z_2} (\Unitary(1),\kap).
    \end{equation*}
    This decomposition 
    immediately implies that, 
    for any cochain
    $\phi_{sc} 
        \in K_\Gam^1(X,(\Spinc(n),\kap_\eps))$,
    there exist cochains
        $\phi_s \in K_\Gam^1(X,(\Spin(n),\id_\eps))$ and
        $\phi_u \in K_\Gam^1(X,(\Unitary(1),\kap_\eps))$
    such that
        $\phi_{sc} = [\phi_s,\phi_u]$.
    It also allows the definition of the map
    \begin{align*}
            \Ad \times q
                : 
            (\Spinc(n),\kap_\eps) 
                &\rightarrow 
            (\SOrth(n),\id_\eps) \times (\Unitary(1),\kap_\eps)\\
                [s,z] &\mapsto (\Ad(s), q(z)).
    \end{align*}
    The next proposition shows that every $\Spink$-structure
    extends to a lifting of a 
    semi-equivariant principal 
    $(\SOrth(n),\id) \times (\Unitary(1),\kap)$-bundle 
    by $\Ad \times q$.

\begin{prop}
\label{prop:QtimesP1}
If $\vphi_0: P \rightarrow Q$ is a $\Spink$-structure,
    then 
    there exists a lifting 
    \begin{equation}
        \vphi: P \rightarrow Q \times_X L
    \end{equation}
    by $\Ad \times q$,
    where $L$ is a $\Gam$-semi-equivariant principal $(\Unitary(1),\kap_\eps)$-bundle.
\begin{proof} 
    Let $\phi \in \TCocycles_\Gam^1(X,(\SOrth(n),\id_\eps))$ be the cocycle for $Q$.
    If $Q$ has a $\Spink$-structure
    there is a cocycle
        $[\phi_s,\phi_u] \in \TCocycles_\Gam^1(X,(\Spinc(n),\kap_\eps))$
    with $\Ad^c([\phi_s,\phi_u]) = \Ad(\phi_s) = \phi$.
    The cocycle $[\phi_s,\phi_u]$ is a lifting by $\Ad \times q$ of $(\phi, \phi_u^2)$.
    It remains to check that $\phi_u^2$ is a cocycle.
First, note that
        $\Ad(\d\phi_s) = \d \circ \Ad(\phi_s) = \d(\phi) = 1$.
    Thus, $\d\phi_s$ takes values in $\ker(\Ad) = \Z_2$, and
    \begin{equation*}
        (\d\phi_s)^\inv(\d\phi_u) 
            \in 
        K_\Gam^2(X,(\Unitary(1),\kap_\eps))
        \subset
        K_\Gam^2(X,(\Spinc(n),\kap_\eps)).
    \end{equation*}
    This cochain is a cocycle as
    \begin{equation*}
        (\d\phi_s)^\inv(\d\phi_u)
            =
        [1, (\d\phi_s)^\inv(\d\phi_u)] 
            =
        [\d\phi_s, \d\phi_u] 
            = 
        \d[\phi_s, \phi_u] 
            = 
        1.
    \end{equation*}
The cochain $\phi_u^2 \in K_\Gam^1(X,(\Unitary(1),\kap_\eps))$ is then a cocycle as
    \begin{align*}
        \d(\phi_u^2)
            =
        (\d\phi_u)^2
            =
        (\d\phi_s)^{-2}(\d\phi_u)^2
            =
        \Big( (\d\phi_s)^\inv(\d\phi_u) \Big)^2
            =
            1.
    \end{align*}
    Therefore, 
    the required bundle $L$ 
    can be constructed from $\phi_u^2$ 
    using
    \cite[Prop.~15]{Kitson.2020.A-Semi-equivariant-Dixmier-Douady-Invariant}.
\end{proof}
\end{prop}

Proposition \ref{prop:QtimesP1} can be refined 
    into a statement about cohomology classes.
    This refinement uses 
    the exact sequences in cohomology
    obtained by applying 
    \cite[Theorem 41]{Kitson.2020.A-Semi-equivariant-Dixmier-Douady-Invariant}
    to the two exact sequences of $\Gam$-groups 
    running diagonally in diagram \eqref{fig:SpincEx}.

\begin{lem}
\label{lem:SOUSeq}
The central exact sequences
    \begin{equation*}
        1 
            \rightarrow  
        (\Z_2,\id_\eps) 
            \rightarrow 
        (\Spin(n),\id_\eps) 
            \os{\Ad}{\rightarrow} 
        (\SOrth(n),\id_\eps)
            \rightarrow 
        1,
    \end{equation*}
    \begin{equation*}
        1 
            \rightarrow  
        (\Z_2,\id_\eps) 
            \rightarrow 
        (\Unitary(1),\kap_\eps) 
            \os{q}{\rightarrow} 
        (\Unitary(1),\kap_\eps) 
            \rightarrow 
        1,
    \end{equation*}
    induce the exact sequences
    \begin{align*}
        H^1_\Gam(X,(\Z_2,\id_\eps))     
          &\os{}{\rightarrow}
        \TCocycles^1_\Gam(X,(\Spin(n),\id_\eps))
          \os{\Ad}{\rightarrow}
            \ldots
        \\
          &\qquad\ldots
          \os{\Ad}{\rightarrow}
        \TCocycles^1_\Gam(X,(\SOrth(n),\id_\eps))
          \os{\Del_{s}}{\rightarrow}
        H^2_\Gam(X,(\Z_2,\id_\eps)),
    \end{align*}
    \begin{align*}
        H^1_\Gam(X,(\Z_2,\id_\eps))     
          &\os{}{\rightarrow}
        H^1_\Gam(X,(\Unitary(1),\kap_\eps))
          \os{q}{\rightarrow}
          \ldots 
        \\
          &\qquad\ldots
          \os{q}{\rightarrow}
        H^1_\Gam(X,(\Unitary(1),\kap_\eps))
          \os{\Del_{u}}{\rightarrow}
        H^2_\Gam(X,(\Z_2,\id_\eps)).
    \end{align*}
\end{lem}

Proposition \ref{prop:QtimesP1} and Lemma \ref{lem:SOUSeq}
    can now be combined to establish an alternative criteria 
    for the existence of a $\Spink$-structure.

\begin{thm}
\label{thm:SOUcond}
A $\Gam$-equivariant principal $\SOrth(n)$-bundle $Q$
    with cocycle $\phi$
    has a $\Spink$-structure
        if and only if
    there exists a cocycle $\psi \in H_\Gam^1(X,(\Unitary(1),\kap_\eps))$
    such that 
      \begin{equation*}
          \Del_s(\phi) = \Del_{u}(\psi) \in H_\Gam^2(X,(\Z_2,\id_\eps)).
      \end{equation*}
\begin{proof} 
Assume that $Q$ has a $\Spink$-structure.
    By Proposition \ref{prop:QtimesP1}, there exists an cocycle
        $[\phi_s,\phi_u] \in \TCocycles_\Gam^1(X,(\Spinc(n),\kap_\eps))$
    such that
        $\phi_u^2$ is a cocycle and 
    \begin{equation*}
        (\Ad \times q)[\phi_s,\phi_u] = (\phi_s, \phi_u^2).
    \end{equation*}
As $[\phi_s,\phi_u]$ is a cocycle,
        $[\d\phi_s,\d\phi_u] = \d[\phi_s,\phi_u] = 1$. 
    This implies that $\d\phi_s = \d\phi_u$.
    Therefore, applying Lemma \ref{lem:SOUSeq}
    to $\phi$ and $\phi_u^2$,
    \begin{equation*}
        \Del_s(\phi) = [\d\phi_s] = [\d\phi_u] = \Del_u(\phi_u^2) \in H_\Gam^2(X,(\Z_2,\id_\eps)).
    \end{equation*}
    Thus, $\psi := \phi_u^2$ is the required cocycle.

    Conversely, suppose there exists a cocycle
        $\psi \in H_\Gam^1(X,(\Unitary(1),\kap_\eps))$ 
    such that
    \begin{equation*}
        \Del_s(\phi) =  \Del_u(\psi) \in H_\Gam^2(X,(\Z_2,\id_\eps)).
    \end{equation*}
    Then, there are a cochains 
        $\phi_s$ with $\Ad(\phi_s)=\phi$, and
        $\phi_u$ with $\phi_u^2=\psi$
    such that
    \begin{equation*}
        [\d\phi_s] = [\d\phi_u] \in K_\Gam^2(X,(\Z_2,\id_\eps)).
    \end{equation*}
    This implies that 
        $\d\phi_s = \d\phi'\d\phi_u = \d(\phi'\phi_u)$ 
    for some $\phi' \in K_\Gam^1(X,(\Z_2,\id_\eps))$.
    Then 
        $\d[\phi_s,\phi'\phi_u] = [\d\phi_s,\d(\phi'\phi_u)] = 1$,
    and
        $\Ad^c[\phi_s,\phi'\phi_u] = \Ad(\phi_s) = \phi$.
    Thus, $[\phi_s,\phi'\phi_u]$ defines a $\Spink$-structure on $Q$.
\end{proof}
\end{thm}


If $X$ is a manifold 
    acted on by 
    a finite group $H$,
    and $V \rightarrow X$
    is a real $H$-equivariant vector bundle
    with cocycle 
    $\phi \in \TCocycles^1_H(X,\SOrth(n))$,
    then the obstruction 
    to the existence of an 
    $H$-equivariant $\Spin$-structure on $V$
    is the 
    second $\Z_2$-valued 
    equivariant 
    Stiefel-Whitney class,
    which can be defined by
    $
    w_2^H(V)
    := 
    \Del_{\Spin}(\phi) \in H^2_H(X,\Z_2)
    $.
    Here $\Del_{\Spin}(\phi)$
    is the connecting map 
    for the exact sequence 
    \begin{equation*}
    \xymatrix{
        H^1_H(X,\Z_2)     
            \ar[r]^(0.425){}
                & 
        \TCocycles^1_H(X,\Spin(n))
            \ar[r]^(0.515){\Ad}
                & 
        \TCocycles^1_H(X,\SOrth(n))
            \ar[r]^(0.55){\Del_\Spin}
                & 
        H^2_H(X,\Z_2), 
             }
    \end{equation*}
    induced 
    by the central exact sequence
    \begin{equation*}
    1 
    \rightarrow 
    \Z_2
    \rightarrow 
    \Spin(n)
    \os{\Ad}{\rightarrow}
    \SO(n)
    \rightarrow 
    1.
    \end{equation*}

If $(\Gam,\eps)$ is the orientifold group
    defined by 
    $\Gam := \Z_2 \times H$
    and
    $\eps(z,h) := z$,
    then $X$ 
    can be made into 
    an orientifold $\til{X}$ for $(\Gam,\eps)$
    by trivially extending 
    its $H$-action
    to the 
    $\Gam$-action $(z,h)x := hx$.
Similarly, 
    the $H$-equivariant vector bundle $V$
    can be made into
    a $\Gam$-equivariant vector bundle $\til{V}$
    by trivially extending 
    its $H$-action
    to the $\Gam$-action 
    $(z,h)v := hv$.
    The cocycle of $\til{V}$
    is an element 
    $\til{\phi} \in 
     \TCocycles^1_\Gam(X,(\SOrth(n),\id_\eps))$.

In this situation,
    the quotient map
    $\pi: \Gam \rightarrow \Gam/\Z_2 \iso H$
    induces a map
    $\pi: X_\Gam^\bul \rightarrow X_H^\bul$
    between the simplicial spaces
    associated to the groups $\Gam$ and $H$.
    Because $\pi$ is a homomorphism and 
    satisfies $\pi(\gam)x = \gam x$,
    it commutes with the face maps
    on these spaces, and
    defines a pulback map $\pi^*$ on cochains.
    The map $\pi^*$ also commutes with 
    the coboundary maps, 
    and provides well-defined 
    extension maps
    \begin{align*}
    \pi^*:
    \TCocycles^p_{H}(X,G)
    &\rightarrow
    \TCocycles^p_\Gam(\til{X},(G,\id_\eps))
    &
    \pi^*:
    H^p_{H}(X,G)
    &\rightarrow
    H^p_\Gam(\til{X},(G,\id_\eps)).
    \end{align*}
One then has the following result.
\begin{prop}
\label{prop:quoComSpin}
If $\til{V} \rightarrow \til{X}$
    is the trivial extension of
    a real $H$-equivariant vector bundle
    $V \rightarrow X$,
    as described above,
    then
\begin{enumerate}
\item 
    the cocycle for $\til{V}$
    is the 
    pullback of 
    the cocycle for $V$ 
    by the quotient map $\pi:\Gam \rightarrow H$,
    \begin{equation*}
    \til{\phi} = \pi^*\phi 
    \in
    H^1_\Gam(\til{X},(\SO(n),\id_\eps)).
    \end{equation*}
\item 
    the second $\Z_2$-valued 
    equivariant Stiefel-Whitney class 
    for $V$ satisfies
    \begin{equation*}
    \pi^* w_2^H(V) 
    = 
    \Del_s(\pi^*\phi)
    \in
    H^2_\Gam(\til{X},(\Z_2,\id_\eps)).
    \end{equation*}
\item 
    $\til{V}$ 
    has a $\Spink$-structure
    if and only if
    \begin{equation*}
    \pi^* w_2^H(V)
    = 
    \Del_{u}(\psi)
    \in
    H^2_\Gam(\til{X},(\Z_2,\id_\eps)),
    \end{equation*}
    for some cocycle 
    $\psi 
        \in 
        H_\Gam^1
            (\til{X},(\Unitary(1),\kap_\eps))$.
\end{enumerate}
Here $\Del_s$ and $\Del_u$ 
    are the connecting maps 
    of Lemma \ref{lem:SOUSeq}.
\begin{proof}
If $\set{s_a}$ is a collection 
    of local sections for $V$,
    then
    \begin{align*}
    \pi(z,h)x &= hx = (z,h)x
    &
    \pi(z,h) s_a(x) &= h s_a(x) = (z,h) s_a(x),
    \end{align*}
    where
    $(z,h) \in \Gam = \Z_2 \times H$,
    $x \in X$.
    Together with
    the property 
    which defines
    the cocycles $\phi$ and $\til{\phi}$
    \cite[Prop.~12]{Kitson.2020.A-Semi-equivariant-Dixmier-Douady-Invariant},
    this implies
    \begin{align*}{}
    s_b(\pi(z,h) x)\phi_{ba}(\pi(z,h),x)
    &=
    \pi(z,h) s_a(x) 
    =
    (z,h) s_a(x) 
    \\&\qquad
    = 
    s_b((z,h) x)\til{\phi}_{ba}((z,h),x)
    \\&\qquad\qquad
    = 
    s_b(\pi(z,h) x)\til{\phi}_{ba}((z,h),x).
    \end{align*}
    Thus, $\pi^*\phi = \til{\phi}$,
    which proves the the first statement.

The second statement follows 
    from the existence 
    of the commutative diagram
    \begin{equation*}
    \xymatrix{
        H^1_H(X,\Z_2)     
            \ar[d]^(0.425){}
            \ar[r]^(0.45){\pi^*}
                & 
        H^1_\Gam(\til{X},(\Z_2,\id_\eps))     
            \ar[d]^(0.425){}
        \\
        \TCocycles^1_H(X,\Spin(n))
            \ar[d]^(0.515){\Ad}
            \ar[r]^(0.45){\pi^*}
                & 
        \TCocycles^1_\Gam(\til{X},(\Spin(n),\id_\eps))
            \ar[d]^(0.515){\Ad}
        \\
        \TCocycles^1_H(X,\SOrth(n))
            \ar[d]^(0.55){\Del_\Spin}
            \ar[r]^(0.45){\pi^*}
                & 
        \TCocycles^1_\Gam(\til{X},(\SOrth(n),\id_\eps))
            \ar[d]^(0.55){\Del_s}
        \\
        H^2_H(X,\Z_2) 
            \ar[r]^(0.45){\pi^*}
                & 
        H^2_\Gam(\til{X},(\Z_2,\id_\eps)).
             }
    \end{equation*}
    To see that 
    the bottom cell of this diagram commutes,
    note that
    if $\psi$ is a lifting of $\phi$,
    then
    $\pi^*\psi$ is a lifting of $\pi^*\phi$.
    The commutation of $\pi^*$ 
    with the coboundary maps
    then implies
    \begin{equation*}
    \pi^* w_2^H(V)
    := 
    \pi^*\Del_\Spin(\phi) 
    =
    \pi^*\d(\psi)
    =
    \d(\pi^*\psi)
    = 
    \Del_s(\pi^*\phi).
    \end{equation*}

The third statement follows 
    from the first and second
    by applying 
    Theorem \ref{thm:SOUcond}.
\end{proof}
\end{prop}


To end this section, 
    two important $\Spink$-structures 
    will be described.
The first of these is 
    the canonical $\Spink$-structure
    associated to a real representation $V$ of
        $(\Z_2,\id) \sdp_{\kap_\eps} \Spinc(n)$.
    When $\dim(V)=8$,
    this $\Spink$-structure is used to construct 
    a canonical reduced orientifold spinor bundle
    over the point orientifold
    for $(\Z_2,\id) \sdp_{\kap_\eps} \Spinc(n)$,
    which, in turn, can be used to construct the
    $8$-fold Bott class over $V$
    for orientifold $K$-theory
    \cite[Example 4.9]{Kitson.2020.Dirac-operators-on-orientifolds}.
The second is a canonical $\Spink$-structure 
    on the $n$-sphere.
    This $\Spink$-structure
    is used to construct 
    a canonical reduced 
    orientifold spinor bundle on $S^{8k}$.
    The reduced 
    orientifold spinor bundle on $S^{8k}$
    can be used to describe
    the compactification 
    of the $8$-fold Bott class 
    over a real representation 
    of $(\Z_2,\id) \sdp_{\kap_\eps} \Spinc(n)$
    \cite[Example 4.11]{Kitson.2020.Dirac-operators-on-orientifolds}.
    .

\begin{lem}[The canonical $\Spink$-structure over a point]
\label{lem:SpinkForV}
Let $V$ be the representation of 
        $(\Z_2,\id) \sdp_{\kap_\eps} \Spinc(n)$
    on $\R^n$ 
    defined by $(\gam,g) \cdot v := \Ad^c(g)v$.
Then 
    \begin{equation*}
        \Ad^c : \Spinc(n) \rightarrow \SO(n) \iso \Fr(V).
    \end{equation*}
    is a $\Spink$-structure
    for the real equivariant vector bundle  
        $V \rightarrow \pt$
    over the point orientifold for 
      $(\Z_2,\id) \sdp_{\kap_\eps} \Spinc(n)$.
\begin{proof}
The group $\Spinc(n)$ forms 
    a principal bundle over a point 
    with the trivial projection $\pi(p) = \pt$,
    and right $\Spinc(n)$ action 
    defined by multiplication.
    The left action of
        $(\Z_2,\id) \sdp_{\kap_\eps} \Spinc(n)$ 
    is taken to be 
    \begin{equation*}
        (\gam,g) \cdot p  := g\kap_{\gam}(p),
    \end{equation*}
    for $\gam \in \Gam$ and $g,p \in \Spinc(n)$.
    The inclusion of the conjugation $\kap$
    is the only difference 
    from the corresponding construction
    in the usual equivariant setting.
\end{proof}
\end{lem}

\begin{lem}[The canonical $\Spink$-structure on the sphere]
\label{lem:SnOSpinc}
The map
    \begin{equation*}
        \Ad^c : \Spinc(n+1) \rightarrow \SO(n+1)
    \end{equation*}
    forms a $\Spink$-structure 
    for the orientifold 
        $S^n \subset \R^{n+1}$
    equipped with the action of
        $(\Z_2,\id) \sdp_{\kap_\eps} \Spinc(n+1)$
    defined by
        $(\gam,g) \cdot v := \Ad^c(g)v$.
\begin{proof}
In what follows, let
        $\gam \in \Z_2$, 
        $g,p \in \Spinc(n+1)$, 
        $h \in \Spinc(n)$, 
        $q \in \SO(n+1)$, 
        $f \in \SO(n)$.
Also, let 
        $\al_1: \SO(n) \rightarrow \SO(n+1)$ and
        $
         \be_1:  \Spinc(n)
               \rightarrow 
               \Spinc(n+1)
        $
    be the maps induced by the inclusion 
        $\Clc_n \rightarrow \Clc_{n+1}$
    defined on standard basis elements by
            $e_k \mapsto e_{k+1}$.
Equip $\Spinc(n+1)$ with the 
    projection, 
    left action, and 
    right $\Spinc(n)$-action
    \begin{align*}
        \pi_{sc}(p)       &:= \Ad^c(p)e_1 &
        (\gam,g) \cdot p  &:= g\kap_{\gam}(p) &
        p \cdot h         &:= p\be_1(h),
    \end{align*}
    respectively.
    Again, the presence 
    of the conjugation action $\kap$ 
    in the left action is the only difference 
    from the corresponding construction
    in the usual equivariant setting
    \cite[p.~5]{Baum-Higson-Schick.2010.A-geometric-description-of-equivariant-K-homology-for-proper-actions}.
    Using the properties of $\kap$, $\Ad^c$ and $\be_1$,
    it is straightforward to check that $\Spinc(n+1)$ forms a
    $(\Gam,\eps) \sdp_{\kap_\eps} \Spinc(n+1)$-semi-equivariant 
    principal $(\Spinc(n),\kap_\eps)$-bundle,
    \begin{align*}
        \pi_{sc}((\gam,g) \cdot p) 
            &= 
        \pi(g(\gam p)) 
        \\&\qquad
            = 
        \Ad^c(g(\gam p))e_1 
            = 
        \Ad^c(g)\Ad^c(\gam p)e_1 
        \\&\qquad\qquad
            = 
        \Ad^c(g)\Ad^c(p)e_1
            = 
        (\gam, g)\pi_{sc}(p),
    \end{align*}
    \begin{align*}
        (\gam,g) \cdot (p \cdot h)
            &=
        (\gam,g) \cdot (p \be_1(h))
        \\&\qquad=
        g( \gam(p \be_1(h)) )
            =
        g (\gam p) (\gam \be_1(h)) 
        \\&\qquad\qquad
            =
        g (\gam p)  \be_1(\gam h) 
        =
        ( (\gam,g)p ) \cdot (\gam h).
    \end{align*}
Next, equip $\SO(n+1)$ with the
    projection, left action, and right $\SO(n)$-action defined by
    \begin{align*}
        \pi_{so}(q)      &:= qe_1 &
        (\gam,g) \cdot q &:= \Ad^c(g)q &
        q \cdot f        &:= q\al_1(f),
    \end{align*}
    respectively.
It can then be checked that $\SO(n+1)$ forms a 
    $(\Z_2,\id) \sdp_{\kap_\eps} \Spinc(n+1)$-equivariant 
    principal $\SO(n)$-bundle,
    \begin{align*}
        \pi_{so}((\gam,g) \cdot q) 
            &= 
        \pi_{so}(\Ad^c(g)q) 
            = 
        \Ad^c(g)qe_1 
            = 
        (\gam, g)\pi(q),
        \\
        (\gam,g) \cdot (q \cdot f)
            &=
        (\gam,g) \cdot (q \al_1(f))
            =
        \Ad(g) q \al_1(f)
            =
        ( (\gam,g) \cdot q ) \cdot f.
    \end{align*}
That $\Ad^c$ is a semi-equivariant lifting can be checked directly
    by verifying compatibility with projections, right actions, and left actions,
    \begin{align*}
        \pi_{sc}(p) 
            &= 
        \Ad^c(p)e_1 
            = 
        \pi_{so} \circ \Ad^c(p),
        \\
        \Ad^c(p \cdot h) 
            &= 
        \Ad^c(p \be_1(h)) 
        \\&\qquad
            = 
        \Ad^c(p)\Ad^c( \be_1(h) )
            = 
        \Ad^c(p)\al_1(\Ad^c(h))
        \\&\qquad\qquad
            = 
        \Ad^c(p) \cdot \Ad^c(h),
        \\
        \Ad^c((\gam,g) \cdot p)
            &=
        \Ad^c(g(\gam p))
        \\&\qquad
            =
        \Ad^c(g)\Ad^c(\gam p)
            =
        \Ad^c(g)\Ad^c(p)
        \\&\qquad\qquad
            =
        (\gam, g) \cdot \Ad^c(p).
    \end{align*}

It remains to check that $\SO(n+1)$ with the given action of
    $(\Z_2,\id) \sdp_{\kap_\eps} \Spinc(n+1)$
    is isomorphic to 
    the equivariant principal $\SO(n)$-bundle $\Fr(S^n)$.
First, identify the tangent space of the $n$-sphere
    with a subbundle of the tangent space to $\R^{n+1}$,
    \begin{equation*}
        TS^n 
            \iso
        \set{
                (v_1,v_2) \in \R^{n+1} \times \R^{n+1} 
                    \mid 
                \Norm{v_1} = \Norm{v_2} = 1, 
                \inn{v_1}{v_2} = 0}
        \subset 
        T\R^{n+1}
    \end{equation*}
The standard action of $\SO(n+1)$ on $\R^{n+1}$
    associates a matrix to each element $q \in \SO(n+1)$,
    which will also be denoted $q$.
    The columns $q_i$ of this matrix determine an orthonormal frame
    \begin{equation*}
        F(q) := \set{(q_1,q_2),\ldots,(q_1,q_{n+1})} \in \Fr_{q_1}(TS^n).
    \end{equation*}
    In this way, $\SO(n+1)$ can be identified with $\Fr(TS^n)$.
    This identification is compatible with projections as
    \begin{equation*}
        \pi_{so}(q) = q e_1 = q_1 = \pi_{TS^n}(F(q)).
    \end{equation*}
    Compatibility with right actions follows from the fact that
    \begin{equation*}
        (q \cdot f)_j 
            = 
        (q \al_1(f))_j 
            =
        \begin{cases}
           q_1 &\txt{ for } j=1 \\
           \sum_{2 \leq i \leq n+1} q_i f_{(i-1)(j-1)} &\txt{ for } j \geq 2.
        \end{cases}
    \end{equation*}
    Finally, the left action on $\Fr(TS^n)$
    can be characterised by observing that
    a vector 
        $(v_1,v) \in TS^n$ 
    is tangent to the curve 
        $(\cos t)v_1 + (\sin t)v$ at $t=0$.
    Acting on this curve by
        $(\gam,g) \in (\Z_2,\id) \sdp_{\kap_\eps} \Spinc(n+1)$
    produces a new curve
        $(\cos t)(\Ad^c(g) v_1) + (\sin t)(\Ad^c(g) v)$
    which has 
        $(\Ad^c(g) v_1,\Ad^c(g) v)$ as its tangent vector at $t=0$.
    Thus,
    \begin{equation*}
        (\gam,g)F(q) = F(\Ad^c(g)q) = F((\gam,g)q),
    \end{equation*}
    and the identification of $\SO(n+1)$ and $\Fr(TS^n)$ is compatible
    with the left actions.
\end{proof}
\end{lem}

\section{Dirac Operators on Orientifolds}
\label{ch:DiracOperators}

In this section, 
    Dirac operators are constructed 
    for orientifolds.
    By applying a semi-equivariant 
    associated bundle construction
    with a Clifford module 
    as the model fibre, 
    it is possible to construct 
    spinor bundles with orientifold actions.
    Both a total spinor bundle, 
    with a right action of $(\Clc_n,\kap_\eps)$,
    and a reduced spinor bundle, 
    with the complexification 
    of an irreducible $\Cl_{8k}$-module 
    as a model fibre, 
    are defined.
    As in the usual setting, 
    the sections of orientifold spinor bundles
    are acted on by 
    sections of a Clifford bundle. 
    This action is compatible with the  
    orientifold action on 
        the spinor bundle and 
        a canonical orientifold
    action on the complex Clifford bundle.
In order to construct a Dirac operator on an orientifold,
    it is neccesary to have a connection 
    which is compatible with
    Clifford multiplication on sections and 
    the orientifold action.
    Such a connection can be constructed 
    using results on semi-equivariant connection 
    forms from \S\ref{sec:SemiConnect}.
After equipping the orientifold spinor bundles 
    with compatible connections, 
    the orientifold Dirac operator and 
    its reduced counterpart will be defined.

\subsection{Orientifold Spinor Bundles}

The model fibre 
    of an orientifold spinor bundle
    is a Clifford module that is
    semi-equivariant 
    with respect to the action 
    of $(\Spinc(n),\kap_\eps)$.
    Such modules can be constructed
    by complexifying $\Cl_n$-modules.
The main $\Cl_n$-modules of interest are $\Cl_n$,
    considered as a module over itself, and
    the irreducible $\Cl_{8k}$-modules.
    Up to equivalence, 
    there is only one 
    irreducible $\Cl_{8k}$-module
    \cite[p.~33]{Lawson-Michelsohn.1989.Spin-geometry}.
    A representative 
    of this equivalence class
    will be denoted by $\Del$.
Denote the complexifications of these,
    equipped with 
    their associated orientifold actions,
    by
    \begin{align*}
        (\Del_c,\kap_\eps) 
        &:= 
        (\Del \ten \C, \id \ten \kap_\eps)
        &
        (\Clc_n,\kap_\eps) 
        &:= 
       (\Cl_n \ten \C,\id \ten \kap_\eps).
    \end{align*}

It is important to note that 
    the complexification
    $\Del \ten \C$ is an irreducible module 
    for $\Clc_{8k}$.
    This is a non-trivial fact that depends 
    on the representation theory 
    of Clifford algberas.

The orientifold spinor bundles can
    now be defined 
    by applying 
    the semi-equivariant 
    associated bundle construction 
    to
    a semi-equivariant principal bundle 
    coming from a $\Spink$-structure
    and a complex Clifford module
    equipped with an orientifold action. 
\begin{defn}\label{defn:DiracBundles}
Let $P \rightarrow \Fr(V)$ be 
    an orientifold-$\Spinc$-structure, and
    define the following orientifold bundles:
\begin{align*}
    &
    \txt{The \emph{orientifold spinor bundle}}
    &
    \tS 
        &:= 
    P 
        \times_{(\Spinc(n),\kap_\eps)} 
    (\Clc_n,\kap_\eps),
    \\
    &
    \txt{The \emph{reduced orientifold spinor bundle}}
    &
    \rS 
    &:= 
    P 
        \times_{(\Spinc(n),\kap_\eps)} 
    (\Del_c,\kap_\eps).
\end{align*}
\end{defn}
Note that if one disregards the orientifold action,
    then an orientifold spinor bundle 
    is a complex spinor bundle in the usual sense.
In the case of the reduced orientifold spinor bundle,
    $\Del_c$ is an irreducible module for $\Clc_{8k}$,
    as mentioned above.
    This implies that, 
    disregarding the orientifold action,
    the reduced orientifold spinor bundle 
    is a reduced complex spinor bundle.

\begin{exam}[The canonical reduced orientifold spinor bundle over a point]
\label{exam:CRSBpt}
Using Lemma \ref{lem:SpinkForV}
    it is possible to construct 
    a $\Spink$-structure $P \rightarrow \Fr(V)$,
    for the adjoint representation $V$ of 
        $(\Z_2,\id) \sdp_{\kap_\eps} \Spinc(n)$.
If $\dim(V) = 8k$, 
    then the irreducible 
    $\Cl_n$-module $\Del$
    can be used to construct 
    a canonical reduced spinor bundle
    $\rS \rightarrow \pt$
    over the point orientifold.
\end{exam}

\begin{exam}[The canonical reduced orientifold spinor bundle over $S^{8k}$]
\label{exam:CRSBS8}
By Lemma \ref{lem:SnOSpinc},
    each sphere $S^n$ has a canonical 
    $(\Z_2,\id) \sdp_{\kap_\eps} \Spinc(n)$-equivariant
    $\Spink$-structure.
If $\dim(V) = 8k$, 
    then the irreducible 
    $\Cl_n$-module $\Del$
    can be used to construct 
    a canonical reduced spinor bundle
    $\rS \rightarrow S^{8k}$
    over the $8$-dimensional sphere.
This construction is an
    adaptation, to the orientifold setting,
    of the Real equivariant spinor bundle
    defined on $S^{8k}$ by Atiyah
    \cite[p.~128]{Atiyah.1968.Bott-periodicity-and-the-index-of-elliptic-operators}.
\end{exam}

The space of sections of the 
    orientifold spinor bundle
    carries an action
    by sections of an 
    \emph{orientifold Clifford bundle} $\Clc(V)$
    called 
    \emph{Clifford multiplication}.
When a $\Spink$-structure 
    $P \rightarrow \Fr(V)$ exists, 
    the orientifold Clifford bundle
    can be expressed as an associated bundle
    \begin{equation*}
        \Clc(V) 
            := 
        P 
            \times^{\Ad^c}_{(\Spinc(n),\kap_\eps)}  
        (\Clc_n,\kap_\eps)
    \end{equation*}
    of $P$,
    and this characterisation 
    can be used to define 
    Clifford multiplication
    on sections of the associated spinor bundle.
    Clifford multiplication on sections is defined
    in terms of the action of $\Clc_n$ 
    on the model fibre.
In order for Clifford multiplication 
    on sections to be well-defined, 
    this fibrewise definition of
    Clifford multiplication 
    must be compatible 
    with the global topology of the base space.
    In the orientifold setting, 
    Clifford multiplication 
    is also required to be compatible 
    with an orientifold action
    on the spinor bundle, and 
    a canonical orientifold action on $\Clc(V)$.
    The $\Spink$-structure used 
    to construct an orientifold spinor bundle
    ensures that 
    both of these requirements are fulfilled.
    Thus, the benefit of working 
    on semi-equivariance and 
    $\Spink$-orientiation
    is finally observed.
In what follows, 
    consider sections 
    of associated bundles to be represented by
    equivariant maps from the principal bundle $P$
    of an underlying $\Spink$-structure 
        $P \rightarrow \Fr(V)$
    into the semi-equivariant fibre, 
    as in Lemma \ref{lem:PSec}.
\begin{prop}\label{prop:CME}
Sections 
        $\vphi \in \Sec(\Clc(V))$ 
    of the orientifold Clifford bundle
    act from the left on the sections 
        $\psi \in \Sec(\tS)$
    of the orientifold spinor bundle
    by
    \begin{equation*}
        (\vphi\psi)(p) = \vphi(p)\psi(p).
    \end{equation*}
    This action is well-defined and satisfies
            $\gam ( \vphi\psi) = (\gam \vphi)(\gam \psi)$.
\begin{proof}
Multiplication is well-defined, 
    as
    \begin{align*}
        (\vphi\psi)(pg) 
            &= \vphi(pg)\psi(pg)
            \\&\qquad
            = (g^\inv \vphi(p)g) (g^\inv\psi(p))
            = g^\inv \vphi(p)\psi(p)
            = g^\inv (\vphi\psi)(p).
    \end{align*}
Compatibility with the orientifold actions 
    is verified using Lemma \ref{lem:PSec},
    \begin{align*}
        (\gam(\vphi\psi))(p)
            &=
        \gam(\vphi\psi)(\gam^\inv p) 
            \\&\qquad
            =
        \gam (\vphi(\gam^\inv p)\psi(\gam^\inv p))
            =
        (\gam \vphi(\gam^\inv p)) (\gam \psi(\gam^\inv p))     
            \\&\qquad\qquad
            =
        (\gam \vphi)(p) (\gam \psi)(p)
            =
        ((\gam \vphi)(\gam \psi))(p).
    \end{align*}
\end{proof}
\end{prop}
Sections of the orientifold Clifford bundle 
    act on sections of the 
    reduced orientifold spinor bundle in the same way.
    One can also check that the 
    Clifford multiplication between sections 
    of the orientifold Clifford bundle
    is well-defined and compatible 
    with the orientifold action.

Because the orientifold spinor bundle 
    has $(\Clc_n,\kap_\eps)$ as its model fibre, 
    it carries a right action 
    by elements of $\Clc_n$.
    This right action 
    is sometimes described as 
    a multigrading 
    \cite[pp.~379-380]{Higson-Roe.2000.Analytic-K-homology}.
\begin{prop}\label{prop:CAE}
    An element $\vphi \in \Clc_n$ acts 
    from the right on sections $\psi \in \Gam(\tS)$
        by
        \begin{equation*}
            (\psi\vphi)(p) = \psi(p)\vphi.
        \end{equation*}
        For $\gam \in \Gam$, this action satisfies
            $\gam (\psi\vphi) = (\gam \psi)(\gam \vphi)$.
\begin{proof}
Consider $\vphi$ as a constant section 
    of the trivial orientifold bundle 
        $P \times^\id_{(G,\tht)} (\Clc_n,\kap_\eps)$.
    The right action is well-defined,
    \begin{equation*}
        (\psi\vphi)(pg) 
            = \psi(pg) \vphi(pg)
            = g^\inv \psi(p)\vphi(p)
            = g^\inv (\psi\vphi)(p).
    \end{equation*}
    It is also compatible with the orientifold actions,
    \begin{align*}
        (\gam(\psi\vphi))(p)
            &=
        \gam(\psi\vphi)(\gam^\inv p) 
        \\&\qquad
            =
        \gam(\psi(\gam^\inv p)\vphi(\gam^\inv p))
            =
        (\gam\psi(\gam^\inv p)) (\gam\vphi(\gam^\inv p))
        \\&\qquad\qquad
            =
        (\gam\psi)(p) (\gam\vphi)(p)
            =
        ((\gam\psi)(\gam\vphi))(p).
    \end{align*}
\end{proof}
\end{prop}
Similar considerations show that there is also 
    a right action of $\Clc_n$ on $\Clc(V)$ 
    which is compatible with their orientifold actions.

\subsection{Connections in Orientifold Spinor Bundles}

In order to define an orientifold Dirac operator,
    a semi-equivariant connection $1$-form 
    is needed
    for the semi-equivariant 
    principal $(\Spinc(n),\kap_\eps)$-bundle
    $P$ of the $\Spink$-structure 
    $P \rightarrow Q$
    underlying the orientifold spinor bundle.
Such a form can be obtained
    by using Proposition \ref{prop:QtimesP1}
    to extend the lifting
        $\vphi: P \rightarrow Q$
    to a lifting 
        $P \rightarrow Q \times_X L$,
    where $L$ is 
    a semi-equivariant principal 
    $(\Unitary(1),\kap_\eps)$-bundle.
A semi-equivariant connection form 
    can then be constructed on 
        $Q \times_X L$, 
    using
    the averaging process of
    Proposition \ref{prop:AvSEConn},
    and lifted to $P$,
    using the relationship 
    between the Lie algebras 
    $\aL{spin^c}(n)$ and $\aL{so}(n) \dsum \aL{u}(1)$.
    In the next proposition, 
    $q$ denotes the square map 
    of Diagram \eqref{fig:SpincEx}.
\begin{prop}
\label{prop:dAdcXq}
The map
    \begin{equation*}
        (\Ad^c \times q)_*
            :
        \aL{spin^c}(n)
            =
        \aL{spin}(n) \dsum \aL{u}(1)
            \rightarrow
        \aL{so}(n) \dsum \aL{u}(1) 
    \end{equation*}
    is an isomorphism, 
    and satisfies
    \begin{equation*}
        (\Ad^c \times q)_* \circ (\id \times \kap_\eps)_*
            =
        (\id \times \kap_\eps)_* \circ (\Ad^c \times q)_*.
    \end{equation*}
\begin{proof}
That 
    $(\Ad^c \times q)_*$ is an isomorphism 
    is a standard result
    \cite[p.~18-20,29]{Friedrich.2000.Dirac-operators-in-Riemannian-geometry}.
    The isomorphism can be written down explicitly
    by making the following identifications
    \begin{enumerate}
    \item 
            $\aL{so}(n)$ 
        can be identified 
        with 
        the real $n \times n$ skew-symmetric matricies.
        A basis for the skew-symmetric matricies
        is defined by
            $\set{E_{ij} \mid 1 \leq i < j \leq n}$
        where 
        $E_{ij}$
        is the $n \times n$ matrix
        with 
            all entries equal to $0$
            except for the
            $(i,j)$th and $(j,i)$th entry,
            which are equal to $1$ and $-1$
            respectively.
    \item
            $\aL{spin}(n)$
        can be identified
        with the linear subspace $\aE^2 \subset \Cl_n$
        spanned by the elements
        $\set{e_ie_j \mid 1 \leq i < j \leq n}$,
        see
        \cite[p.~18]{Friedrich.2000.Dirac-operators-in-Riemannian-geometry}.
    \item
        $\aL{u}(1)$ 
        can be identified with $\R$.
    \end{enumerate}
With these identifications,
        $(\Ad^c \times q)_*$
    is the map
    \begin{align*}
        (\Ad^c \times q)_*
            :
        \aL{spin}(n) \dsum \aL{u}(1)
            &\rightarrow
        \aL{so}(n) \dsum \aL{u}(1) 
    \\
        (e_ie_j,t) &\mapsto (2E_{ij},2t),
    \end{align*}
    see 
    \cite[pp.~19-20,29]{Friedrich.2000.Dirac-operators-in-Riemannian-geometry}.
    Also, the $\Gam$-actions
    on 
        $\aL{spin}(n) \dsum \aL{u}(1)$ 
        and
        $\aL{so}(n) \dsum \aL{u}(1)$ 
    are 
    \begin{align*}
        (\id \dsum \kap_\eps)_*: 
            \aL{spin}(n) \dsum \aL{u}(1) 
                &\rightarrow
            \aL{spin}(n) \dsum \aL{u}(1)
    \\
        (e_ie_j,t) &\mapsto (e_ie_j,\iot_\eps(t))
    \end{align*}
    \begin{align*}
        (\id \dsum \kap_\eps)_*: 
            \aL{so}(n) \dsum \aL{u}(1) 
                &\rightarrow
            \aL{so}(n) \dsum \aL{u}(1)
    \\
        (E_{ij},t) &\mapsto (E_{ij},\iot_\eps(t)),
    \end{align*}
    where
    $\iot_\eps: \R \rightarrow \R$
    is the involutive action 
    induced by $\iot: t \mapsto -t \in \R$.
Examining these maps,
    it is clear that
        $
            (\Ad^c \times q)_* 
                \circ 
            (\id \times \kap_\eps)_*
                =
            (\id \times \kap_\eps)_* 
                \circ 
            (\Ad^c \times q)_*
        $.
\end{proof}
\end{prop}

\begin{prop}
\label{prop:SEconnform}
Let $\vphi_Q: P \rightarrow Q$ be a $\Spink$-structure. 
    The semi-equivariant principal bundle $P$ 
    carries a $\Gam$-semi-equivariant connection $1$-form.
\begin{proof}
By Proposition \ref{prop:QtimesP1}, 
    there exists 
        a lifting
        \begin{equation*}
            \vphi_Q \times \vphi_L
                : 
            P 
                \rightarrow 
            Q \times_X L
        \end{equation*}
        by $\Ad^c \times q$,
        where $L$ is 
        a semi-equivariant principal 
        $(\Unitary(1),\kap_\eps)$-bundle.
The equivariant principal bundle $Q$ has 
    an equivariant connection $1$-form 
        $\omg_Q: TQ \rightarrow \aL{so}(n)$
        determined by an equivariant metric.
The semi-equivariant principal bundle $L$
    has a semi-equivariant connection $1$-form
    $\omg_L : TL \rightarrow \aL{u}(1)$
    constructed 
    by 
    applying Proposition \ref{prop:AvSEConn}
    to any choice of connection $1$-form for $L$.
Together, these two connection $1$-forms define
    a semi-equivariant connection $1$-form
    \begin{equation*}
        \omg_Q \dsum \omg_L
            : 
        T(Q \times_X L) 
            \rightarrow 
        \aL{so}(n) \dsum \aL{u}(1).
    \end{equation*}
Using the $\Spink$-structure $\vphi$
    and Proposition \ref{prop:dAdcXq},
    the connection $1$-form
        $\omg_Q \dsum \omg_L$ 
    can be lifted to a connection $1$-form
    \begin{align*}
        \omg
            : 
        TP 
            &\rightarrow 
        \aL{spin^c}(n) 
    \\
        v 
            &\mapsto 
        (\Ad^c \times q)_*^\inv 
            \circ 
        (\omg_Q \dsum \omg_L) 
            \circ 
        (\vphi_Q \times \vphi_L)_*(v).
    \end{align*}
    The semi-equivariance of $\omg$ 
    follows from the 
    semi-equivariance of $\omg_Q \dsum \omg_L$,
    and the 
    equivariance of
    $(\vphi_Q \times \vphi_L)_*$ and
    $(\Ad^c \times q)_*$.
\end{proof}      
\end{prop}

The next proposition shows that 
    the connection $1$-form
    constructed by Proposition \ref{prop:SEconnform}
    defines 
    a covariant derivative
    on the orientifold spinor bundle
    that is equivariant 
    with respect to the action of $\Gam$.
    In this proposition,
    sections will be considered
    as maps 
        $\psi: P \rightarrow \Clc_n$
    satisfying 
        $\psi(gp) = g^\inv \psi(p)$,
    and will be acted on 
    by the $\Gam$-action defined 
    in Lemma \ref{lem:PSec}.
    From the point of view 
    of the exterior covariant derivative, 
    these maps are order-zero tensorial forms
        $\psi \in \Lam^0(P,\Clc_n)$.  
    For the details of tensorial forms
    and exterior covariant derivatives,
    see 
    \cite[\S B.3-4]{Friedrich.2000.Dirac-operators-in-Riemannian-geometry}
    \cite[\S II.5]{Kobayashi-Nomizu.1963.Foundations-of-differential-geometry.1}.

\begin{prop}
\label{prop:SemiConnEqui}
Let $\vphi: P \rightarrow Q$ be a $\Spink$-structure. 
    The semi-equivariant connection $1$-form $\omg$,
    defined on $P$ by Proposition \ref{prop:SEconnform},
    determines 
    an exterior covariant derivative
    \begin{equation*}
        d^\omg
            : 
        \Lam^0(P,\Clc_n)  
            \rightarrow 
        \Lam^1(P,\Clc_n)
    \end{equation*}
    that satisfies the condition
    \begin{equation*}
        d^\omg(
            \kap_{\eps(\gam)} 
                \circ 
            \psi 
                \circ 
            \eta_{\gam^\inv}
        ) 
            = 
        \kap_{\eps(\gam)} 
            \circ 
        d^\omg\psi 
            \circ 
        (\eta_{\gam^\inv})_*,
    \end{equation*}
    where 
        $\psi \in \Lam^0(P,\Clc_n)$,
        $\eta$ is the $\Gam$-action on $P$, and 
        $\kap_\eps$ is the conjugation action on $\Clc_n$.
\begin{proof}
The vertical projection
    associated to 
    the connection form $\omg$ 
    is defined by
    \begin{equation*}
       \pi_V|_p 
            := 
        (R^p)_* \circ \omg
            : 
        TP_p \rightarrow TP_p.
    \end{equation*}
    Therefore,
    the exterior covariant derivative
    can be written as
    \begin{equation}
    \label{eq:ExCoD}
        d^\omg\psi(v) 
            = 
        d\psi \circ \pi_H(v) 
            = 
        d\psi(v) -  d\psi \circ \pi_V(v) 
            =
        d\psi(v) - d\psi \circ (R^p)_* \circ \omg(v),
    \end{equation}
    where 
        $v \in TP_p$,
        $\psi \in \Lam^0(P,\Clc_n)$, and
        $\pi_H$ is the horizontal projection.
The first term of the decomposition \eqref{eq:ExCoD} 
    is equivariant,
    as the properties of the exterior derivative
    imply that
    \begin{equation*}
        d(
            \kap_{\eps(\gam)} 
            \circ 
            \psi 
            \circ 
            \eta_{\gam^\inv}
         ) 
            = 
        \kap_{\eps(\gam)} 
        \circ 
        d\psi 
        \circ 
        (\eta_{\gam^\inv})_*.
    \end{equation*}
The semi-equivariance of $P$
    implies the identity
    $(\eta_{\gam})_* \circ (R^p)_* 
        = (R^{\gam p})_* \circ (\tht_{\gam})_*$.
    Together with the 
    the semi-equivariance of $\omg$,
    this implies that
    \begin{align*}
        d(
            \kap_{\eps(\gam)} 
                \circ 
            \psi 
                \circ 
            \eta_{\gam^\inv}
        ) 
            \circ 
        (R^p)_* 
            \circ 
        \omg
            &=
        \kap_{\eps(\gam)}
            \circ 
        d\psi
            \circ 
        (\eta_{\gam^\inv})_* 
            \circ 
        (R^p)_* \circ \omg 
    \\
            &=
        \kap_{\eps(\gam)} 
            \circ 
        d\psi
            \circ 
        (R^{\gam^\inv p})_* 
            \circ 
        (\tht_{\gam^\inv})_* 
            \circ 
        \omg 
    \\
            &=
        \kap_{\eps(\gam)} 
            \circ 
        d\psi
            \circ 
        (R^{\gam^\inv p})_* 
            \circ  
        \omg 
            \circ 
        (\eta_{\gam^\inv})_*.
    \end{align*}
    Therefore, 
    the second term of the decomposition \eqref{eq:ExCoD} 
    is also equivariant.
\end{proof}
\end{prop}
Proposition \ref{prop:SemiConnEqui}
    applies equally well to the reduced
    orientifold spinor bundle 
    if $\Clc_n$ is replaced with $\Del_c$.

As in the non-equivariant case,
    the exterior covariant derivative is
    also equivariant with respect
    to the right action of $\Clc_n$ on 
    the orientifold spinor bundle.

\begin{prop}
\label{prop:RightSemiConnEqui}
Let $\vphi: P \rightarrow Q$ be a $\Spink$-structure. 
    The semi-equivariant connection $1$-form $\omg$,
    defined on $P$ by Proposition \ref{prop:SEconnform},
    determines 
    an exterior covariant derivative
    \begin{equation*}
        d^\omg
            : 
        \Lam^0(P,\Clc_n)  
            \rightarrow 
        \Lam^1(P,\Clc_n)
    \end{equation*}
    that satisfies
    \begin{equation*}
        d^\omg(\psi\vphi) 
            = 
        d^\omg(\psi)\vphi, 
    \end{equation*}
    for 
        $\psi \in \Lam^0(P,\Clc_n)$ and 
        $\vphi \in \Clc_n$.
\end{prop}

\subsection{Dirac Operators on Orientifolds}
\label{sec:DOOOF}

At this stage, 
    all of the preliminary constructions 
    have been completed.
    It is now possible 
    to construct the orientifold Dirac operator
    and reduced orientifold Dirac operator.

\begin{defn}
\label{defn:ODO}
    Let $\conn^L$ denote the connection associated 
    to a $\Spink$-structure 
        $P \rightarrow \Fr(TM)$ 
    by Proposition \ref{prop:SEconnform},
    and $\mu$ denote Clifford multiplication 
    by sections of 
        $T^*M \iso TM \subset \Clc(TM)$.
    Define the 
    \emph{orientifold Dirac operator}
    and
    \emph{reduced orientifold Dirac operator},
    respectively,
    by
\begin{align*}
    \tD
        &:= 
    \mu \circ \conn^L 
        : 
    \Gam(\tS) 
        \rightarrow 
    \Gam(T^*M \ten \tS) 
        \rightarrow 
    \Gam(\tS),
    \\
    \rD
        &:= 
    \mu \circ \conn^L
        : 
    \Gam(\rS) 
        \rightarrow 
    \Gam(T^*M \ten \rS) 
        \rightarrow 
    \Gam(\rS).
\end{align*}
\end{defn}

The 
    orientifold Dirac operator and 
    reduced orientifold Dirac operator
    are complex Dirac operators, 
    in the usual sense. 
    However, 
    they are equivariant 
    with respect to the orientifold actions
    on their spinor bundles.    
    Thus, when 
        $\eps: \Gam \rightarrow \Z_2$ 
    is non-trivial, they have anti-linear symmetries.

\begin{prop}\label{prop:ODLequi}
The orientifold Dirac operator is equivariant with respect to 
    the left action of $\Gam$ on sections of $\tS$,
    \begin{equation*}
        \tD(\gam \psi) = \gam \tD(\psi).
    \end{equation*}
\begin{proof}
This follows from Propositions \ref{prop:CME} and \ref{prop:SemiConnEqui}.
\end{proof}
\end{prop}
The same arguments show that 
    the reduced orientifold spinor bundle
    is also $\Gam$-equivariant.
In addition to $\Gam$-equivariance, 
    the orientifold Dirac operator is equivariant 
    with respect to 
    the right action of $(\Clc_n,\kap_\eps)$ 
    on the orientifold spinor bundle.
\begin{prop}\label{prop:ODRequi}
The orientifold Dirac operator is equivariant with respect to 
    the right action of $\Clc_n$ on sections of $\tS$,
    \begin{equation*}
       \tD(\psi\vphi) = \tD(\psi)\vphi.
    \end{equation*}
\begin{proof}
This follows from Propositions \ref{prop:CAE} and \ref{prop:RightSemiConnEqui}. 
\end{proof}
\end{prop}
Note, in particular, 
    that left and right equivariance together
    imply that the index of $\tD$ consists of 
    vector spaces which are both Clifford modules
    and orientifold representations of $(\Gam,\eps)$.

The main aim of this paper is now complete,
    and the following theorem has been proved.
\begin{thm}
\label{thm:RealDirac}
Let $X$ be an orientifold 
    with orientifold group $(\Gam,\eps)$.
\begin{enumerate}
\item 
    If $W^{(\Gam,\eps)}_3(X)=0$, then 
    $X$ carries 
    an orientifold Dirac operator.
\item 
    If $W^{(\Gam,\eps)}_3(X)=0$
    and $\dim(X) = 8$, 
    then 
    $X$ carries 
    a reduced orientifold Dirac operator.
\end{enumerate}
In particular,
    if $X$ is an 
    $8k$-dimensional Real manifold
    and $W^{(\Z_2,\id)}_3(X)=0$, 
    then $X$ carries a 
    reduced Real Dirac operator.
\end{thm}

\section{Related Work and Applications}
\label{sec:RWAPP}

To put the construction 
    of the orientifold Dirac operator 
    in context,
    it is worth breifly
    recalling the position that
    $\Spin$-structures and 
    $\Spin$-Dirac operators 
    occupy in 
    the $K$-theory of real vector bundles.
The centrality of $\Spin$-structures
    and the $\Spin$-Dirac operator
    in the $K$-theory of real vector bundles
    stems from their
    role in Atiyah's index theoretic proof 
    of the Thom isomorphism theorem
    \cite{Atiyah.1968.Bott-periodicity-and-the-index-of-elliptic-operators},
    which provides an isomorphism 
    $KO(X) \rightarrow KO(V)$
    between
    the $KO$-theory 
    of a manifold
    and any $8k$-dimensional $\Spin$-oriented  
    real vector bundle $V \rightarrow X$.
The proof proceeds
    by compactifying the fibres of $V$ into 
    a family of $8k$-dimensional spheres, 
    each of which is equipped 
    with a canonical reduced Dirac operator.
    The families index map 
    $KO(V) \rightarrow KO(X)$
    associated to 
    this family of Dirac operators 
    is then shown to provide 
    an inverse to the Thom map.
This method of proof naturally accomodates 
    additional symmetries,
    provided that
    the appropriate analogue 
    of $\Spin$-structure and Dirac operator
    can be determined.
    In this way,
    Atiyah was able to prove
    the equivariant Thom isomorphism
    \cite{Atiyah.1968.Bott-periodicity-and-the-index-of-elliptic-operators}.
The $\Spink$-structures and 
    orientifold Dirac operators
    defined in this paper
    play an exactly analogous role
    in the $K$-theory of orientifold bundles.
    These allow
    Atiyah's argument to be extended,
    providing a proof of
    the corresponding Thom isomorphism theorem
    for orientifold $K$-theory
    \cite[Theorem 4.30]{Kitson.2020.Dirac-operators-on-orientifolds}.
    They also provide a basis
    for the definition of 
    geometric orientifold $K$-homology
    \cite[Chapter 6]{Kitson.2020.Dirac-operators-on-orientifolds}.

Orientation conditions in $K$-theory, 
    such as $\Spink$-orientibility,
    are closely related to the topic 
    of twisted $K$-theory.
    The twisted $K$-theory 
    of Real topological groupoids
    has been studied by
    Moutuou using 
    a \v{C}ech cohomology for Real groupoids
    \cite{Moutuou.2013.On-groupoids-with-involution-and-their-cohomology}.
    In a more algebraic context,
    Karoubi and Weibel
    have studied 
    an equivariant 
    twisted $K$-theory that
    includes $KR$-theory as a special case
    \cite{Karoubi-Weibel.2017.Twisted-K-theory-real-A-bundles-and-Grothendieck-Witt-groups}.
    Another approach
    by
    Hekmati et al.
    \cite{Hekmati-Murray-Szabo-Vozzo.2016.Real-bundle-gerbes-orientifolds-and-twisted-KR-homology},
    motivated by applications 
    to orientifold string theories,
    studies $KR$-orientiation and twisting
    using a Real sheaf cohomology theory.
    The work of Freed and Moore 
    on topological insulators
    \cite{Freed-Moore.2013.Twisted-equivariant-matter}
    also treats twistings of $K$-theory
    in the presence of symmetry.
    
Orientifold string theory and 
    the classification 
    of topological insulators
    are 
    two areas in which 
    the constructions of this paper 
    have potential applications.
The connection between the present investigation and
    string theory begins with the classification of
    D-brane charges using $K$-theory, as described in 
    \cite{Minasian-Moore.1997.K-theory-and-Ramond-Ramond-charge, Witten.1998.D-branes-and-K-theory}.
    Results in index theory allow one to pass 
    from $K$-theory to an analytic $K$-homology theory
    in which 
    classes are represented 
    by elliptic operators.
    Each class in this $K$-homology theory
    may be represented 
    by a Dirac operator
    that has been twisted with a vector bundle.
    By replacing these Dirac operators
    with classes formed from
    the $\Spinc$-structures and vector bundles
    used to construct them,
    it is possible to define a 
    $K$-homology theory 
    in entirely geometric terms
    \cite{Baum-Douglas.1982.K-homology-and-index-theory,
          Baum-Douglas.1982.Index-theory-bordism-and-K-homology,
          Baum-Higson-Schick.2007.On-the-equivalence-of-geometric-and-analytic-k-homology}.
    This characterisation of D-brane charge
    is of interest, 
    as the geometric data 
    comprising such a $K$-homology 
    class has physical interpretations
    \cite{Baum.2010.K-homology-and-D-branes}
    \cite[\S 4]{Szabo.2002.D-branes-Tachyons-and-K-homology}.
Three types of orientifold string theories
    are listed in 
    \cite[p.~26-27]{Witten.1998.D-branes-and-K-theory},
    along with the corresponding $K$-theories 
    that classifying the 
    associated D-brane charges.
  In the first of these, 
    D-brane charges 
    are classified by $KR$-theory.
    The $\Spink$-structure and orientifold
    Dirac operators constructed
    in this paper provide the 
    ingredients neccesary to 
    generalise the above discussion
    to $KR$-theory and 
    the $K$-theory of orientifold bundles.
  The two other possibilities 
    listed in 
    \cite[p.~26-27]{Witten.1998.D-branes-and-K-theory}
    involve $K$-theory with sign-choice.
    This $K$-theory has been studied 
    by Doran et al.
    \cite{Doran-Mendez-Diez-Rosenberg.2014.T-duality-for-orientifolds-and-twisted-KR-theory} 
    using methods from non-commutative geometry.
    The 
    $K$-theory with sign-choice is a subgroup
    of the $K$-theory orientifold bundles. 
    Many of the constructions 
    discussed in the present paper
    could be modified to incorporate
    sign-choice structures.

In recent years, 
    there there has been much interest in 
    the classification of topological insulators.
    These classification attempts lead naturally 
    to the consideration 
    of topological invariants 
    which respect anti-linear symmetries
    \cite{Kane-Mele.2005.Z2-topological-order-and-the-quantum-spin-hall-effect,
    Fu-Kane.2006.Time-reversal-polarization-and-a-Z2-adiabatic-spin-pump,
    Fu-Kane.2007.Topological-insulators-with-inversion-symmetry,
    Fu-Kane-Mele.2007.Topological-insulators-in-three-dimensions,
    Qi-Hughes-Zhang.2008.Topological-field-theory-of-time-reversal-invariant-insulators}.
    Contact with Clifford algebras and $K$-theory
    has been made
    through the work of Kitaev
    \cite{Kitaev.2009.Periodic-table-for-topological-insulators-and-superconductors}.
    Another framework 
    for studying topological insulators,
    using twisted $K$-theories, 
    has been described by Freed and Moore
    \cite{Freed-Moore.2013.Twisted-equivariant-matter}.
    The $K$-theory of orientifold bundles
    is a primary example within their framework.
    Thus, it appears that there is potential
    for index invariants derived from
    the orientifold Dirac operator
    to be applied to the classification 
    of topological insulators.

\vfil

\appendix
\section{Semi-equivariant Constructions}

\subsection{Semi-equivariance and Associated Bundles}

The construction 
    of associated bundles 
    from semi-equivairant principal bundles
    differs slightly 
    from the corresponding 
    equivariant construction.
    When forming an equivariant vector bundle 
    from an equivariant principal bundle,
    the only requirement 
    on the model fibre is that it carries
    carries an action of the structure group $G$.
    However, 
    when forming a vector bundle 
    from a semi-equivariant principal bundle,
    it is neccesary to use a model fibre
    that carries both 
    an action of the structure group $G$ and 
    an action of the equivariance group $\Gam$.
    As on the semi-equivariant 
    principal bundle,
    these two actions 
    are required to combine 
    into an action 
    of the 
    semi-direct product group $\Gam \sdp_\tht G$.
    Although the action 
    of the equivariance group $G$
    on the model fibre 
    is required to be linear,
    the action 
    of the equivariance group $\Gam$ is not.
    This makes it possible
    to construct associated bundles
    with $\Gam$-actions that are not linear.
    In particular,
    it is possible to construct
    complex vector bundles 
    equipped with linear/anti-linear actions
    as semi-equivariant 
    associated bundles.

\begin{defn}
Let $P$ be 
    a $\Gam$-semi-equivariant 
    principal $(G,\tht)$-bundle.
A \emph{semi-equivariant fibre} for $P$ 
    is a vector space $V$ 
    equipped with 
    a linear action of $G$
    and an action of $\Gam$ 
    by diffeomorphisms,
    such that 
    \begin{equation*}
        \gam(gv) = (\gam g)(\gam v).
    \end{equation*}
\end{defn}

\begin{defn}
\label{defn:SEAB}
Let 
    $P$ be 
    a $\Gam$-semi-equivariant 
    principal $(G,\tht)$-bundle, and 
    $V$ be a semi-equivariant fibre for $P$.
The \emph{semi-equivariant associated bundle} 
    is the vector bundle
    \begin{equation*}
        P \times_{(G,\tht)} V := P \times V / \sim
    \end{equation*}
    where $(p,v) \sim (pg^\inv,gv)$.
    This bundle carries an action of $\Gam$ defined by
    \begin{equation*}
        \gam(p,v) := (\gam p,\gam v).
    \end{equation*}
\end{defn}
Note that the $\Gam$-action on $P \times_{(G,\tht)} V$ 
    is well-defined because
    \begin{align*}
        \gam[pg^\inv, gv] 
            &= 
        [\gam (pg^\inv), \gam (gv)] 
        \\&\qquad
            = 
        [(\gam p)(\gam g)^\inv, (\gam g)(\gam v)]
            = 
        [\gam p,\gam v] 
            = 
        \gam[p,v].
    \end{align*}

Sections of associated bundles 
    are often represented as equivariant maps
    from the principal bundle 
    into the model fibre.
    It is sometimes useful to express 
    the action of $\Gam$ on a section 
    in this way.
\begin{lem}
\label{lem:PSec}
    Sections of $P \times_{(G,\tht)} V$ 
    are in bijective correspondence with maps
        $\psi: P \rightarrow V$ 
    such that 
        $\psi(pg) = g^\inv \psi(p)$.
    The $\Gam$-action on sections 
    of $P \times_{(G,\tht)} V$
    corresponds to the 
    $\Gam$-action
    \begin{equation*}
        (\gam\psi)(p) = \gam \psi(\gam^\inv p)
    \end{equation*}
    on these maps.
\begin{proof}
    A map 
        $\psi: P \rightarrow V$ 
    with 
        $\psi(pg) = g^\inv \psi(p)$
    corresponds to 
    the section 
    of $P \times_{(G,\tht)} V$ 
    defined by 
        $s(p) := [p,\psi(p)]$.
    The $\Gam$-action on such a section is
    \begin{align*}
        (\gam s)(p) 
            &:= 
        \gam s(\gam^\inv p) 
        \\&\qquad
            = 
        \gam [\gam^\inv p,\psi(\gam^\inv p)] 
            =
         [\gam\gam^\inv p,\gam \psi(\gam^\inv p)] 
            =
         [p,\gam \psi(\gam^\inv p)].
    \end{align*}
    Thus, the corresponding map on $P$ is $p \mapsto \gam\psi(\gam^\inv p)$.
\end{proof}
\end{lem}

\subsection{Semi-equivariant Connections}
\label{sec:SemiConnect}

In the smooth non-equivariant setting, 
    a connection for a principal $G$-bundle $P$ 
    can be expressed 
    as a $\aL{g}$-valued $1$-form 
    on the tangent space $TP$,
    where $\aL{g}$ is the Lie algebra 
    of the structure group $G$
    \cite[Chapter 2]{Kobayashi-Nomizu.1963.Foundations-of-differential-geometry.1},
    \cite[Appendix B]{Friedrich.2000.Dirac-operators-in-Riemannian-geometry}.
A $\Gam$-semi-equivariant 
    $(G,\tht)$-principal bundle 
    \cite[\S 2]{Kitson.2020.A-Semi-equivariant-Dixmier-Douady-Invariant}
    has a $\Gam$-group $(G,\tht)$ 
    as its structure group.
    The differentials $(\tht_\gam)_*$ 
    of the $\Gam$-action on $G$
    form a $\Gam$-action 
    on the Lie algebra $\aL{g}$.
    A connection in 
    a semi-equivariant principal bundle 
    must be compatible
    with this action 
    if it is to produce 
    an equivariant connection 
    in an associated bundle.
    The definition of 
    a semi-equivariant connection $1$-form 
    is given below, 
    along with an averaging proceedure 
    that can be used to construct 
    semi-equivariant connections.
    In what follows, let 
        $R_g(p) = R^p(g) := pg$ 
    denote the multiplication maps 
    associated to the right action 
    on a principal $G$-bundle $P$.
    Also, 
    let $R_g(h):=hg$
    denote 
    the right action of $G$ on itself.
    Note that $(R^p)_*(A_e)$ 
    defines 
    the vector field 
    induced on $P$ by an element $A \in \aL{g}$,
    and the adjoint map on $\aL{g}$ 
    may be expressed as $\Ad_{g^\inv} = (R_g)_*$.
\begin{defn}
Let 
    $(P,\eta)$ be a 
    smooth $\Gam$-semi-equivariant 
    principal $(G,\tht)$-bundle 
    with $\Gam$-action $\eta$,
    and
    let $\aL{g}$ be the Lie alegebra of $G$.
A \emph{$\Gam$-semi-equivariant connection $1$-form} on $P$ 
    is a Lie algebra valued $1$-form 
    \begin{equation*}
        \omg: TP \rightarrow \aL{g}
    \end{equation*}
    such that 
    for all 
        $\gam \in \Gam$, 
        $g \in G$, 
        $A \in \aL{g}$, and 
        $p \in P$,
    \begin{align*}
        \omg \circ (R^p)_*(A_e) 
        &= 
        A 
    &
        \omg \circ (R_g)_* 
        &= 
        (R_g)_* \circ \omg  
    &
        \omg \circ (\eta_\gam)_* 
        &= 
        (\tht_\gam)_* \circ \omg.
    \end{align*}
\end{defn}

When $\Gam$ is finite, a semi-equivariant connection can be constructed from
    a given connection by a twisted averaging procedure.

\begin{prop}\label{prop:AvSEConn}
Let $\Gam$ be a finite Lie group, and 
    suppose that $P$ is 
    a smooth $\Gam$-semi-equivariant 
    principal $(G,\tht)$-bundle 
    with $\Gam$-action $\eta$.
    If 
        $\omg: TP \rightarrow \aL{g}$
    is a connection form on $P$, then
    \begin{equation*}
        \omg_\Gam 
        := 
        \sum_{\gam \in \Gam} 
            (\tht_\gam)_* 
            \circ 
            \omg 
            \circ 
            (\eta_{\gam^\inv})_*
    \end{equation*}
    is a $\Gam$-semi-equivariant connection on $P$.
\begin{proof}
First note that, as
    $\tht$ is an automorphism and 
    $P$ is semi-equivariant,
    identities are induced between the differentials of the various actions.
    For $\gam \in \Gam$, $g,h \in G$, and $p \in P$
    \begin{alignat*}{3}
        \gam(hg) 
        &= 
        (\gam h)(\gam g) 
            &&\implies 
            \hspace{1em} 
        (\tht_\gam)_* \circ (R_g)_* 
        = 
        (R_{\gam g})_* \circ (\tht_\gam)_*  
    \\
        \gam(pg)  
        &= 
        (\gam p)(\gam g)  
            &&\implies 
        \begin{cases}
            (\eta_\gam)_* \circ (R_g)_* 
            = 
            (R_{\gam g})_* \circ (\eta_\gam)_*  
        \\
            (\eta_\gam)_* \circ (R^p)_* 
            = 
            (R^{\gam p})_* \circ (\tht_\gam)_*.
        \end{cases}
    \end{alignat*}
To check that $\omg_\Gam$ is a connection, 
    first observe that the condition
        $\omg_\Gam \circ (R^p)_*(A_e) = A$
    holds,
    \begin{align*}
        (\tht_\gam)_* \circ \omg 
        \circ 
        (\eta_{\gam^\inv})_* \circ (R^p)_*(A_e)
        &= 
        (\tht_\gam)_* \circ \omg 
        \circ 
        (R^{\gam^\inv p})_* \circ (\tht_{\gam^\inv})_*(A_e) 
    \\
        &= 
        (\tht_\gam) \circ \omg 
        \circ 
        (R^{\gam^\inv p})_* ( (\tht_{\gam^\inv})_*(A)_e ) 
    \\
        &= 
        (\tht_\gam)_* \circ (\tht_{\gam^\inv})_*(A) 
    \\
        &= 
        A.
    \end{align*}
    The condition 
        $\omg_\Gam \circ (R_g)_* 
            = 
        (R_g)_* \circ \omg_\Gam$ 
    also holds, as
    \begin{align*}
        (\tht_\gam)_* \circ \omg 
        \circ 
        (\eta_{\gam^\inv})_* \circ (R_g)_*
        &= 
        (\tht_\gam)_* \circ \omg 
        \circ 
        (R_{\gam^\inv g})_* \circ (\eta_{\gam^\inv})_*  
    \\
        &= 
        (\tht_\gam)_* \circ (R_{\gam^\inv g})_* 
        \circ 
        \omg \circ (\eta_{\gam^\inv})_*        
    \\ 
        &= 
        (R_g)_* \circ (\tht_\gam)_* 
        \circ 
        \omg \circ (\eta_{\gam^\inv})_*.
    \end{align*}
Finally, semi-equivariance holds, as 
    \begin{align*}
        \omg_\Gam \circ (\eta_{\gam})_* 
        &= 
        (
         \sum_{\gam_1 \in \Gam} 
            (\tht_{\gam_1})_* \circ \omg 
            \circ 
            (\eta_{\gam_1^\inv})_*
        ) 
        \circ 
        (\eta_{\gam})_* 
    \\
        &= 
        \sum_{\gam_1 \in \Gam} 
            (\tht_{\gam_1})_* \circ \omg 
            \circ 
            (\eta_{\gam_1^\inv \gam})_* 
    \\
        &= 
        \sum_{\gam_2 \in \Gam} 
            (\tht_{\gam \gam_2^\inv})_* \circ \omg 
            \circ 
            (\eta_{\gam_2})_* 
    \\
        &= 
        (\tht_{\gam})_* 
            \circ 
        (
         \sum_{\gam_2 \in \Gam}  
            (\tht_{\gam_2^\inv})_* \circ \omg 
            \circ 
            (\eta_{\gam_2})_*) 
    \\
        &= 
        (\tht_{\gam})_* \circ \omg_\Gam.
    \end{align*}
\end{proof}
\end{prop}